\RequirePackage{etoolbox}
\csdef{input@path}{{style/}{graphics/}}
\documentclass[seceqn,MSNbibl,number,citesort,dvips]{arxbj}
\usepackage{mathrsfs,upgreek,newsym2e}
\usepackage{graphicx}

\aid{0}
\volume{21}
\issue{1}
\pubyear{2015}
\firstpage{144}
\lastpage{175}
\doi{10.3150/13-BEJ563} 

\makeatletter

\newcommand{\rright}{\right}
\newcommand{\lleft}{\left}
\newcommand{\rrvert}{\vert}
\newcommand{\llvert}{\vert}
\newcommand{\N}{\mathbb{N}}%

\newcommand{\R}{\mathbb{R}}%

\newcommand{\E}{\mathbb{E}}%
\newcommand{\Prob}{\mathbb{P}}%

\newcommand{\eqref}[1]{(\ref{#1})}
\newcommand{\supp}{\operatorname{supp}}
\newcommand{\Var}{\operatorname{Var}}
\renewcommand{\backslash}{\setminus}
\renewcommand{\pi}{\uppi}

\newcommand{\binom}[2]{\pmatrix{#1\cr #2}}
\newcommand{\bigtimes}{\mathop{\!\mbox{\parbox[c][9pt][b]{18pt}{\fontsize{18}{18}\selectfont{$\times$}}}\!\!\!\!}}
\newtheorem{Theorem}{Theorem}
\newtheorem{Corollary}{Corollary}
\newtheorem{Lemma}{Lemma}
\newproclaim{Assumption}{Assumption}
\newremark{Example}{Example}
\newremark{Remark}{Remark}

\def\sfrac#1#2{#1/#2}
\def\vfrac#1#2{(#1)/#2}
\def\afrac#1#2{#1/(#2)}

\makeatother

\begin{document}
\begin{frontmatter}

\title{Confidence bands for multivariate and time dependent inverse
regression models}
\runtitle{Confidence bands for multivariate regression models}

\begin{aug}
\author[A]{\inits{K.}\fnms{Katharina} \snm{Proksch}\thanksref{e1}\ead[label=e1,mark]{katharina.proksch@rub.de}},
\author[A]{\inits{N.}\fnms{Nicolai} \snm{Bissantz}\thanksref{e2}\ead[label=e2,mark]{nicolai.bissantz@rub.de}} \and
\author[A]{\inits{H.}\fnms{Holger} \snm{Dette}\thanksref{e3}\ead[label=e3,mark]{holger.dette@rub.de}}
\address[A]{Fakult\"at f\"ur Mathematik, Ruhr-Universit\"at Bochum,
44780 Bochum, Germany.\\ \printead{e1,e2,e3}}
\end{aug}

\received{\smonth{6} \syear{2012}}
\revised{\smonth{5} \syear{2013}}

%
\begin{abstract}
Uniform asymptotic confidence bands for a
multivariate regression function
in an inverse regression model
with a convolution-type operator are constructed.
The results are derived using strong approximation methods and a
limit theorem for the supremum of a stationary Gaussian field over an
increasing system of sets.
As a particular application, asymptotic confidence bands for a
time dependent regression function $f_t(x)$ ($x\in\R^{d}, t\in\R$)
in a convolution-type inverse
regression model are obtained.
Finally, we demonstrate the practical feasibility of our proposed
methods in a simulation study and an
application to the estimation of the luminosity profile of the
elliptical galaxy NGC5017.
To the best knowledge of the authors, the results presented in this
paper are the first which provide
uniform confidence bands for multivariate nonparametric function
estimation in inverse problems.
\end{abstract}

%
\begin{keyword}
\kwd{confidence bands}
\kwd{deconvolution}
\kwd{inverse problems}
\kwd{multivariate regression}
\kwd{nonparametric regression}
\kwd{rates of convergence}
\kwd{time dependent regression function}
\kwd{uniform convergence}
\end{keyword}

\end{frontmatter}

\section{Introduction}
\setcounter{equation}{0}

\subsection{Inverse regression models} \label{sec1}
In many applications, it is impossible to observe a certain quantity of
interest because only indirect observations are available for
statistical inference. Problems of this type are called inverse
problems and arise in many fields such as medical imaging, physics and
biology. Mathematically the connection between the quantity of interest
and the observable one can often be expressed in terms of a linear
operator equation.
Well-known examples are Positron Emission Tomography, which involves
the Radon Transform (Cavalier \cite{cavalier2000}),
the heat equation (Mair and Ruymgaart \cite{mairuy1996}), the
Laplace Transform (Saitoh \cite{saitoh1997}) and the
reconstruction of astronomical and biological images from telescopic
and microscopic imaging devices, which is closely connected to
convolution-type operators (Adorf \cite{adorf1995},
Bertero \textit{et al.} \cite{berbocdes2009}).

Inverse problems have been studied intensively in a deterministic
framework and in mathematical physics. See, for example, Engl \textit{et al.} \cite{enghanneu1996} for an overview of existing
methods in numerical
analysis of inverse problems or Saitoh \cite
{saitoh1997} for techniques based on
reproducing kernel Hilbert spaces.
Recently, the investigation of inverse problems has also become of
importance from a statistical point of view. Here, a particularly interesting
and active field of research is the construction of statistical
inference methods such as hypothesis tests or confidence regions.

In this paper, we are interested in the convolution type inverse
regression model
%
\begin{equation}
\label{mod1} Y = (f*\psi) (x) + \varepsilon,
\end{equation}
where $\varepsilon$ is a random error, the operation $*$ denotes
convolution, $\psi$ is a given square integrable function and the
object of interest is the function $f$ itself.
An important and interesting application of the inverse regression
model (\ref{mod1}) is the recovery of images from imaging
devices such as astronomical telescopes or fluorescence microscopes in
biology. In these cases, the observed,
uncorrected image is always at least slightly blurry due to the
physical characteristics of the propagation of
light at surfaces of mirrors and lenses in the telescope. In this
application, the
variable $x$ represents the pixel of a CCD and
we can only observe a blurred version of the true image modeled by the
function $f$.
In the corresponding mathematical model,
the observed image is (at least approximately) a convolution of the
real image with the so-called
point spread function $\psi$, that is, an inverse problem with
convolution operator.

The inference problem regarding the function $f$ is called inverse
problem with stochastic noise. In recent years, the problem of estimating
the regression function $f$ has become an important field of research,
where the main focus is on a one dimensional predictor. Several authors
propose Bayesian methods (Bertero \textit{et al.} \cite
{berbocdes2009}, Kaipio and Somersalo \cite{kaisom2005})
and construct estimators using tools from nonparametric curve
estimation (Mair and Ruymgaart \cite{mairuy1996}, Cavalier \cite{cavalier2008}, Bissantz \textit{et al.}
\cite{bishohmunruy2007}). Further inference methods, in particular the
construction of confidence intervals and confidence bands, are much
less developed. Birke \textit{et al.} \cite{birbishol2010} have
constructed uniform
confidence bands for the
function $f$ with a one-dimensional predictor.

The present work is motivated by the fact that in many applications one
has to deal with an at least two-dimensional
predictor. A typical example is image reconstruction since a picture is
a two-dimensional object.
Also in addition to the spatial dimensions, the data often show a
dynamical behavior, thus repeated measurements
at different times can be used to extend the statistical inference. For example,
in astrophysics spectra of different objects like supernovae or
variable stars show changes in time on observable timescales.
In this case, the function $f$ depends on a further parameter, say
$f_t$ and the reconstruction problem refers to a multivariate function
even if the predictor is univariate.

The purpose of the present paper is the investigation of asymptotic
properties of estimators for the function $f$
in model \eqref{mod1} with a multivariate predictor. {In particular,
we present a result on the weak convergence
of the sup-norm of an appropriately centered estimate, which can be
used to construct
asymptotic confidence bands for the regression function $f$.}
In contrast to other authors (e.g., Cavalier and Tsybakov \cite
{cavtsy2002}), we do not
assume that the function $\psi$ in model \eqref{mod1}
is periodic, because in the reconstruction of astronomical or
biological images from telescopes or microscopic imaging devices
this assumption is often unrealistic.

\subsection{Confidence bands}\label{sec2}
In a pioneering work, Bickel and Rosenblatt \cite{bicros1973a}
extended results of Smirnov \cite{smirnov1950} for a
histogram estimate and constructed confidence bands
for a density function of independent, identically distributed (i.i.d.)
observations. Their method is based on the asymptotic distribution of
the supremum of a centered kernel density estimator. Since then, their
method has been further developed both in the context of density and
regression estimation. For density estimation, Neumann \cite{neumann1998}
derived bootstrap confidence bands, and Gin{\'{e}} and Nickl \cite
{ginnic2010} derived
adaptive asymptotic bands over generic sets. In a regression context,
asymptotic confidence bands were constructed by Eubank and Speckman \cite{eubspe1993} for
the Nadaraya--Watson estimator and by Xia \cite
{xia1998} for a local
polynomial estimator. Bootstrap confidence bands for nonparametric
regression were proposed by Hall \cite{hall1993},
Neumann and Polzehl \cite{neupol1998} and by
Claeskens and {van Keilegom} \cite{clakei2003}. For the statistical
inverse problem of
deconvolution density estimation, Bissantz \textit{et al.} \cite
{bisdumholmun2007} constructed
asymptotic and bootstrap confidence bands, where Lounici and Nickl \cite{Lounic2011} obtained
non-asymptotic confidence bands by using concentration inequalities.
Recently, Birke \textit{et al.} \cite{birbishol2010} provided
uniform asymptotic and
bootstrap confidence bands for a spectral cut-off estimator in the
one-dimensional indirect regression model with convolution operator.

All these results are limited to the estimation of univariate densities
and regression functions, and are not applicable in cases where the
quantity of interest depends on a multivariate predictor. In such
cases, to the best knowledge of the authors, confidence bands are not
available. One reason for this gap is that a well-established way to
construct asymptotic uniform confidence bands, which uses a pioneering
result of
Bickel and Rosenblatt \cite{bicros1973a} as the standard tool,
cannot be extended in a
straightforward manner to the multivariate case. There are substantial
differences between the multivariate and one-dimensional case, and for
multivariate inverse problems the mathematical construction of
confidence bands requires different and/or extended methodology.

In the present paper, we will consider the problem of constructing
confidence bands for the regression function in an inverse regression
model with a convolution-type operator with a multivariate predictor.
The estimators and assumptions for our asymptotic theory are presented
in Section~\ref{sec3}, while Section~\ref{sec4} contains the main
results of the paper. In Section~\ref{sec5}, we consider the special
case of time dependent regression functions with a univariate
predictor, which originally motivated our investigations. In Section~\ref{sec:sim},
the finite-sample properties of the proposed asymptotic confidence
bands are illustrated by means of a small simulation study and an
application to HST data is discussed.
The arguments of Sections~\ref{sec6} and \ref{sec7}, which contain
all technical details of the proofs, are based on results by
Piterbarg \cite{piterbarg1996} who
provided a limit theorem
for the supremum
\[
\sup_{t\in T_n} X(t)
\]
of a stationary Gaussian field $\{X(t) | t \in\mathbb{R}^d\} $,
where $\{T_n\subset\R^d\}_{n\in\N}$
is an increasing system of sets such that $\lambda^d(T_n)\rightarrow
\infty$ as $n\rightarrow\infty$.
This result generalized the multivariate extension in
Bickel and Rosenblatt \cite{bicros1973b}, who provided a
limit theorem for the supremum
$\sup_{t\in[0,T]^d} X(t)$,
as $T\rightarrow\infty$.

\section{Notation and assumptions}\label{sec3}
\setcounter{equation}{0}

\subsection{Model and notations} \label{prel}
Suppose that $(2n+1)^d$ observations $(x_{\mathbf{k}},Y_{\mathbf
{k}}),\mathbf{k}=(k_1,\ldots, k_d)\in G_{\mathbf{n}} := \{-n,\ldots
,n\}^d$ from the model
%
\begin{eqnarray}
\label{model} Y_{\mathbf{k}}=g(x_{\mathbf{k}})+\varepsilon_{\mathbf{k}}:=
(f*\psi ) (x_{\mathbf{k}})+\varepsilon_{\mathbf{k}},
\end{eqnarray}
are available, where the function $f\dvtx \R^d\rightarrow\R$ is unknown,
$\psi\dvtx \R^d\rightarrow\R$ is a known function and $
g:= f*\psi$ denotes the convolution of $f$ and $\psi$, that is
%
\begin{eqnarray}
\label{conv} g(x):= (f*\psi) (x):=\int_{\mathbb{R}^d}f(s)\psi(x-s)
\,\mathrm{d}s.
\end{eqnarray}
The basic assumptions that guarantee the existence of the integral
(\ref{conv}) and also assure $g\in L^2(\R^d)$ is that $f\in L^2(\R
^d)$ and $\psi\in L^1(\R^d)\cap L^2(\R^d)$, which will be assumed
throughout this paper.
In model \eqref{model}, the predictors
$x_{\mathbf{k}}:= \mathbf{k}\cdot\frac{1}{na_n}$ are equally spaced
fixed design points on a $d$-dimensional grid, with a
sequence $(a_n)_{n\in\N}$ satisfying
\begin{eqnarray*}
na_n\rightarrow\infty\quad \mbox{and}\quad  a_n\searrow0
\qquad \mbox{for } n\rightarrow\infty.
\end{eqnarray*}
The noise terms $\{\varepsilon_{\mathbf{k}} | \mathbf{k}\in
G_{\mathbf{n}}\}$ are a field of centered i.i.d. random variables with
variance $\sigma^2:=\mathbb{E}\varepsilon_{\mathbf{k}}^2>0$ and
existing fourth moments.
As a consequence of the convolution theorem and the formula for Fourier
inversion, we obtain the representation
%
\begin{eqnarray}
\label{repf} f(x)=\frac{1}{(2\pi)^d}\int_{\R^d}
\frac{ \mathcal{F}g(\xi
)}{\mathcal{F}\psi(\xi)}\exp\bigl(\mathrm{i}\xi^Tx\bigr) \,\mathrm{d}\xi.
\end{eqnarray}
An estimator for the regression function $f$ can now easily be obtained
by replacing the unknown quantity $\mathcal{F}g=\mathcal{F} (f*\psi
)$ by an estimator $ \mathcal{F}\hat g$. The random fluctuations in
the estimator $ \mathcal{F}\hat g$ cause instability of the ratio
$\frac{\mathcal{F}\hat g(\xi)}{\mathcal{F}\psi(\xi)}$ if at least
one of the components of $\xi$ is large.
As a\vspace*{1pt} consequence, the problem at hand is ill-posed and requires
regularization. We address this issue by excluding large values of $\xi
_j$ for any $j=1,\ldots,d$ from the domain of integration, that is, we
multiply the integrand in (\ref{repf})
with a sequence of Fourier transforms $\mathcal{F} \eta(h\cdot)$ of
smooth functions with compact support $[-h^{-1},h^{-1}]^d$. Here
$h=h_n$ is a regularization parameter which corresponds to a bandwidth
in nonparametric curve estimation and satisfies $h\rightarrow0$ if
$n\rightarrow\infty$.
For the exact properties of the function $\eta$, we refer to
Assumption \ref{assA} below.

An estimator $\hat f_n$ for the function $f$ in model (\ref{model}) is
now easily obtained as
%
\begin{eqnarray}
\label{est} \hat f_n(x)=\frac{1}{(2\pi)^d}\int
_{\R^d}\frac{ \mathcal{F}\hat
g(\xi)}{\mathcal{F}\psi(\xi)}\exp\bigl(\mathrm{i}\xi^Tx\bigr)
\mathcal{F}\eta(h\xi ) \,\mathrm{d}\xi,
\end{eqnarray}
where
\begin{eqnarray*}
\mathcal{F}\hat g(\xi)=\frac{1}{(2\pi)^{\sfrac{d}{2}}n^da_n^d}\sum_{\mathbf{k}\in G_n}Y_{\mathbf{k}}
\exp\bigl(-\mathrm{i}\xi^Tx_{\mathbf{k}}\bigr)
\end{eqnarray*}
is the empirical analogue of the Fourier transform of $g$.
Note that with the definition of the kernel
%
\begin{eqnarray}
\label{kernKn} K_n(x)=\frac{1}{(2\pi)^{\sfrac{d}{2}}}\int_{\R^d}
\frac{\mathcal
{F}\eta(\xi)}{\mathcal{F}\psi(\sfrac{\xi}{h})}\exp\bigl(\mathrm{i}\xi^Tx\bigr) \,\mathrm{d}\xi,
\end{eqnarray}
the estimator (\ref{est}) has the following representation
%
\begin{eqnarray}
\label{esti} \hat f_n(x)=\frac{1}{(2\pi)^{d}n^da_n^d h^d}\sum
_{\mathbf{k}\in
G_\mathbf{n}}Y_{\mathbf{k}}K_n \biggl((x-x_{\mathbf{k}})
\frac
{1}{h} \biggr).
\end{eqnarray}
Note that the kernel $K_n$ can be expressed as a Fourier transform as follows
\begin{eqnarray*}
K_n=\overline{\mathcal{F} \biggl(\frac{\mathcal{F}\eta}{\mathcal
{F}\psi(\sfrac{\cdot}{h})} \biggr)}.
\end{eqnarray*}
Also note that the kernel $K_n$ is a so-called deconvolution kernel. It
is the analogue of a kernel in classical nonparametric kernel
estimation with the difference that it depends on $n$ via the bandwidth
$h$ in a rather complicated manner. For this reason, we use the
notation $K_n$ instead of $K_h$
which corresponds to $K_h(\cdot/h) =\frac{1}{h^d}K(\cdot/h)$. Asymptotically,
this kernel can be replaced by  its limit $K$, see Assumption \ref{assB},
Remark~\ref{remark2} and Example~\ref{Ex1} in the following
discussion.

The first step of the proof of our main result (see Theorem~\ref{GWS}
in Section~\ref{sec4}) will consist of a uniform approximation of
$\hat f_n(x)-\E\hat f_n(x)$ by an appropriate stationary Gaussian field.
In the second step, we apply results of Piterbarg
\cite{piterbarg1996} and Bickel and Rosenblatt \cite{bicros1973b}
to obtain the desired uniform convergence for the approximation process
of the first step. Finally, these results are used to construct uniform
confidence regions for $\E\hat f_n(x)$.
Our approach is then based on undersmoothing: the choice of
sufficiently small bandwidths assures the same limiting behaviour of
$\hat f_n(x)-\E\hat f_n(x)$ and $\hat f_n(x)-f(x)$. This avoids the
estimation of higher order derivatives, which often turns out to be
difficult in applications.
Thus, the limit theorem obtained in the second step will also provide
uniform confidence regions for the function $f$ itself. Whereas
undersmoothing implies that the rate-optimal bandwidth cannot be used,
there has also been some theoretical justification why this choice of
the regularization parameter is useful for constructing confidence
intervals (see Hall \cite{hall1992b}).

\subsection{Assumptions} \label{subsecass}
We now introduce the necessary assumptions which are required for the
proofs of our
main results in Section~\ref{sec4}. The first assumption refers to
the type of (inverse) deconvolution problem describing
the shape of the kernel function $\eta$ in the spectral domain.

\renewcommand{\theAssumption}{\Alph{Assumption}}
\begin{Assumption}\label{assA}
Let $\mathcal{F}\eta$ denote the Fourier
transform of a function $\eta$ such that
\begin{enumerate}[A3.]
\item[A1.] $\supp (\mathcal{F}\eta )\subset[-1,1]^d$.
\item[A2.] $\mathcal{F}\eta\in\mathscr{D}(\R^d)=\{f\dvtx \R
^d\rightarrow\R| f\in C^{\infty}
(\R^d), \supp(f) \subset\R^d \mbox{ compact}\}$.
\item[A3.] There exists a constant $D>0$, such that $\mathcal{F}\eta
(\xi)=1$ for all $\xi\in[-D,D]^d$ and
$|\mathcal{F}\eta(\xi)|\leq1$ for all $\xi\in\R^d$.
\end{enumerate}
\end{Assumption}

\begin{Remark}
\begin{enumerate}[3.]
\item[1.] The decay of the tails of the kernel $K_n$ is given in terms
of the smoothness of the integrand in (\ref{kernKn}). The choice of a
smooth regularizing
function $\mathcal{F}\eta$ has the advantage that the smoothness of
$1/\mathcal{F}\psi$ carries over to $\mathcal{F}\eta(h\cdot)
/\mathcal{F}\psi$.
\item[2.] Functions like $\mathcal{F}\eta$ are called bump
functions. Their existence follows from the $C^{\infty}$ Urysohn lemma
(see, e.g., Folland \cite{folland1984}, Lemma~8.18).
\item[3.] Note that $\mathscr{D}(\R^d)\subset\mathscr{S}(\R^d)$,
where $\mathscr{S}(\R^d)$ denotes the Schwartz space of smooth and
rapidly decreasing functions. Since $\mathcal{F}\dvtx \mathscr{S}(\R
^d)\rightarrow\mathscr{S}(\R^d)$ is a bijection (see, e.g.,
Folland \cite{folland1984}, Corollary~8.28)
we know that $\eta\in\mathscr{S}(\R^d)$
as well.
\item[4.] For the sake of transparency, we state the conditions and
results with the same regularization parameter $h$ for each direction.
In practical applications, this might not be the best strategy. The
results presented in Sections~\ref{sec4} and \ref{sec5} also hold for
different sequences of bandwidths $h_1,\ldots,h_d$ as long as the
system of rectangles
$\{[0,h_1^{-1}]\times\cdots\times[0,h_d^{-1}] | n\in\N\}$ is a blowing
up system of sets in the sense of Definition~14.1 in Piterbarg \cite{piterbarg1996}. This is the case if the assumption
\begin{eqnarray*}
\sum_{p=1}^{d} \Biggl(\prod
_{j=1,j\neq p}^{d}\frac{1}{h_j} \Biggr)\leq
L_1 \cdot \Biggl(\prod_{j=1}^d
\frac{1}{h_j} \Biggr)^{L_2},
\end{eqnarray*}
is satisfied for a constant $L_1$ that only depends on $d$ and a
constant $L_2<1$.
This condition is not a restriction in our setting because it holds
whenever $h_j\cdot n^{\gamma_j}\rightarrow C_j$ for constants $C_j,
\gamma_j>0$, $ j=1,\ldots,d$.
\end{enumerate}
\end{Remark}

In general, two kinds of convolution problems are distinguished in the
literature, because the decay of the Fourier transform of the
convolution function $\psi$ determines the degree of ill-posedness. In
the case of an exponentially decreasing Fourier transform $\mathcal
{F}\psi$ the problem is called severely ill-posed. In the present
paper, the class of moderately ill-posed problems is considered, where
the Fourier
transform of the convolution function decays at a polynomial rate (the
precise condition will be specified in Assumption \ref{assB} below).
Throughout this paper
\begin{eqnarray*}
\mathcal{W}^m \bigl(\mathbb{R}^d\bigr) =\bigl\{f\in
L^2\bigl(\R^d\bigr) | \partial^{(\alpha
)}f\in
L^2\bigl(\R^d\bigr) \mbox{ exists } \forall\alpha\in
\N^d, |\alpha |\leq m \bigr\},
\end{eqnarray*}
denotes the Sobolev space of order $m \in\mathbb{N}$, where $\partial
^{(\alpha)}f$ is the
weak derivative of $f$ of order $\alpha$. In the subsequent
discussion, we will also make use of
the Sobolev space for general $m>0$, which is defined by
\begin{eqnarray*}
\mathcal{W}^m \bigl(\mathbb{R}^d\bigr) =\bigl\{f\in
L^2\bigl(\R^d\bigr) | \bigl(1+|\xi|^2
\bigr)^{\sfrac
{m}{2}}\mathcal{F}f\in L^2\bigl(\R^d\bigr)
\bigr\}.
\end{eqnarray*}

\begin{Assumption}\label{assB} We assume the existence of a function $\Psi\dvtx \R
^d\rightarrow\R$ such that the kernel $K=\overline{\mathcal{F}
(\Psi\cdot\mathcal{F}\eta )}$ satisfies
\begin{enumerate}[B3.]
\item[B1.] $K\neq0$ and
there exist constants $\beta> d/2$, $M\in\N$, indices $0<\mu_1<\mu
_2<\cdots<\mu_M$ and $L^2$-functions
$f_1,\ldots,f_{M-1},f_{M}\dvtx \R^d\rightarrow\R$ with the property
\[
\xi^{\alpha}f_p\in\mathcal{W}^m\bigl(
\R^d\bigr) \qquad (p=1,\ldots,M-1)
\]
for all multi-indices $\alpha\in\{0,\dots,d\}^d, |\alpha|\leq d$
and all
$m>\frac{d+|\alpha|}{2}$, such that
%
\begin{eqnarray}
\label{expan} h^{\beta}K_n(x)-K(x)=\sum
_{p=1}^{M-1}h^{\mu_p}\overline{\mathcal
{F}f_p}(x)+h^{\mu_M}\overline{\mathcal{F}f_{n,M}}(x),
\end{eqnarray}
where $f_M$ may depend on $n$, that is, $f_M=f_{M,n}$ and $\|f_{M,n}\|
_{L^1(\R^d)}=\mathrm{O}(1)$.

\item[B2.] $\xi^{\alpha}\Psi\cdot\mathcal{F}\eta, \xi^{\alpha
}\frac{h^{\beta}}{\mathcal{F}\psi(\sfrac{\cdot}{h})}\cdot\mathcal
{F}\eta\in\mathcal{W}^m(\R^d)$ for some $m>\frac{d+|\alpha|}{2}$.

\item[B3.] $\log(n)\cdot h^{\mu_M}(a_n^{-\sfrac{d}{2}}h^{-\sfrac
{d}{2}})\cdot\|f_{M}\|_{L^1(\R^d)}=\mathrm{o}(1)$ and $h^{\mu_{1}}(\log(n))^2=\mathrm{o}(1)$.
\end{enumerate}
\end{Assumption}

\begin{Remark}\label{remark2}Assumption B1 implies $h^{\beta
}K_n\rightarrow K$ in $L^2(\R^d)$
and also specifies the order of this convergence. It can be understood
as follows. If the convergence of
the difference $h^{\beta}K_n-K$ is fast enough, that is,
%
\begin{eqnarray}
\label{case1} \log(n)\cdot h^{\mu_1}(a_nh)^{-\sfrac{d}{2}}=\mathrm{o}(1)
\end{eqnarray}
we have $M=1$.
On the other hand, in some relevant situations (see Example~\ref{Ex1}(ii) below)
the rate of convergence $h^{\mu_1}$ is given by $h^2$ for each $d$ and
\eqref{case1} cannot
hold for $d\geq4$. Here, the expansion (\ref{expan}) provides a
structure, such that our main
results remain correct although the rate of convergence is not very
fast. We can decompose the
difference $h^{\beta}K_n-K$ in two parts, where one part depends on
$n$ only through the factors
$h^{\mu_p}$ and the other part converges sufficiently fast (in some
cases this term vanishes completely).
\end{Remark}

\begin{Example}\label{Ex1} This example illustrates the construction
of the functions in the representation \eqref{expan}.
\begin{enumerate}[(ii)]
\item[(i)] Let $d=2$ and $\psi(x)=\frac{1}{4}\exp(-|x_1|)\exp
(-|x_2|)$, $x=(x_1,x_2)^T$, $\xi=(\xi_1,\xi_2)^T$. Then we have
$\frac{h^4}{\mathcal{F}\psi(\xi/h)}=2\pi (h^4+h^2(\xi
_1^2+\xi_2^2)+\xi_1^2\xi_2^2 )$,
which implies $\beta=4, M=3$ and
\begin{eqnarray*}
h^4\cdot K_n(x)&=&\int_{\R^2}
\bigl(h^4+h^2\bigl(\xi_1^2+
\xi_2^2\bigr)+\xi _1^2
\xi_2^2 \bigr)\mathcal{F}\eta(\xi)\exp
\bigl(\mathrm{i}x^T\xi\bigr) \,\mathrm{d}\xi,
\\
K(x)&=&\int_{\R^2}\mathcal{F}\eta(\xi)
\xi_1^2\xi_2^2\exp
\bigl(\mathrm{i}x^T\xi \bigr) \,\mathrm{d}\xi.
\end{eqnarray*}
With the definitions $f_1(\xi)=2\pi(\xi_1^2+\xi_2^2)\mathcal
{F}\eta(\xi)$, $f_2(\xi)=2\pi\mathcal{F}\eta(\xi)$ and
$f_{n,3}\equiv0$ we obtain
\begin{eqnarray*}
h^4\cdot K_n(x)-K(x)=h^2\cdot\overline{
\mathcal{F}f_1(\xi )}+h^4\cdot\overline{
\mathcal{F}f_2(\xi)}.
\end{eqnarray*}
In this example, the condition $\log(n)h^2/\sqrt{a_n^dh^d}=\mathrm{o}(1)$ is
satisfied. However, the following results are valid if the weaker
condition of a decomposition of the form (\ref{expan}) holds.
Furthermore, since the factors of $\mathcal{F} \eta$ in $f_1$ and
$f_2$ are polynomials, we have $\mathcal{F}f_j(\xi)\in\mathscr
{S}(\R^d)$, which implies $\xi^{\alpha}f_j\in\mathcal{W}^m(\R^d)$
for all $\alpha$ and all $m\in\N$.
\item[(ii)] If $|x|=\sqrt{x_1^2+\cdots+x_d^2}$ and $\psi
(x)=2^{-\vfrac{d+1}{2}}\mathrm{e}^{-|x|}$ we have
\begin{eqnarray*}
\mathcal{F}\psi(\xi)=\frac{1}{\sqrt{2\pi}}\Gamma \biggl(\frac
{d+1}{2} \biggr)
\frac{1}{(1+|\xi|^2)^{\vfrac{d+1}{2}}},
\end{eqnarray*}
(see Folland \cite{folland1984}, Exercise 13).
If $d$ is odd we use the identity
\begin{eqnarray*}
\bigl(h^{2}+|\xi|^2 \bigr)^{\vfrac{d+1}{2}}=\sum
_{j=0}^{\vfrac
{d+1}{2}}\binom{\frac{d+1}{2}\vspace*{2pt}}
{j}h^{2j}|\xi|^{\vfrac{d+1}{2}-2j},
\end{eqnarray*}
and an expansion of the form \eqref{expan} is obvious from the
definition of $K_n$ in \eqref{kernKn}. If the dimension $d$ is even
the situation is more complicated. Consider for example the case $d=4$, where
\begin{eqnarray*}
\frac{h^{5}}{\mathcal{F}\psi(\sfrac{\xi}{h})} \rightarrow\frac{\sqrt{2\pi}}{\Gamma (\sfrac{5}{2}
)}|\xi|^{5} =
\frac{\sqrt{2\pi}}{\Gamma (\sfrac{5}{2} )}\sqrt{ \bigl(\xi_1^2+
\xi_2^2+\xi_3^2+
\xi_4^2 \bigr)^{5}}\qquad  \mbox{as } n\rightarrow
\infty.
\end{eqnarray*}
It follows that the constant $\beta$ and the functions $\Psi, K_n$
and $K$ from Assumption \ref{assB} are given by $\beta=d+1=5$, $\Psi(\xi)
=\frac{\sqrt{2\pi}}{\Gamma (\sfrac{5}{2} )}|\xi|^{5}$ and
\begin{eqnarray*}
h^{\beta}K_n(x) &=&\frac{1}{(2\pi)^{2}}\int
_{\R^d}\frac{\sqrt{2\pi}}{\Gamma
(\sfrac{5}{2} )} \bigl(h^{2}+|
\xi|^2 \bigr)^{\vfrac{d+1}{2}} \mathcal{F}\eta(\xi)\exp\bigl(\mathrm{i}
\xi^Tx\bigr) \,\mathrm{d}\xi,
\\
K(x)&=&\frac{1}{(2\pi)^{2}}\int_{\R^d}\frac{\sqrt{2\pi}}{\Gamma
(\sfrac{5}{2} )}|
\xi|^{d+1} \mathcal{F}\eta(\xi)\exp\bigl(\mathrm{i}\xi^Tx\bigr) \,\mathrm{d}\xi,
\end{eqnarray*}
respectively.
In order to show that Assumption B1 holds in this case we use
Taylor's theorem and obtain
\begin{eqnarray*}
\frac{h^{5}}{\mathcal{F}\psi(\sfrac{\xi}{h})}-\Psi(\xi)= \frac{\sqrt{2\pi}}{\Gamma (\vfrac{d+1}{2} )} \biggl( h^2\cdot
\frac{5}{2}\cdot|\xi|^3+h^4\cdot\frac{5}{2}
\cdot\frac
{3}{2}\cdot\bigl(|\xi|^2+\lambda_dh^2
\bigr)^{\sfrac{1}{2}} \biggr),
\end{eqnarray*}
for some constant $\lambda_d\in[0,1)$.
Recalling the definition of $K_n$ in \eqref{kernKn} this gives
\begin{eqnarray*}
\bigl(h^{\beta}K_n-K\bigr) (x)&= h^2 \overline{
\mathcal{F}f_1(\xi)}+h^4 \overline{\mathcal
{F}f_{2,n}(\xi)},
\end{eqnarray*}
where the functions $f_1$ and $f_{2,n}$ are defined by
\begin{eqnarray*}
f_1(\xi) & =&\frac{1}{(2\pi)^{\sfrac{3}{2}}\Gamma (\sfrac
{5}{2} )}\cdot|\xi|^3\cdot
\frac{5}{2}\cdot\mathcal{F}\eta (\xi),
\\
f_{2,n}(x) & =& \frac{1}{(2\pi)^{\sfrac{3}{2}}\Gamma (\sfrac{5}{2}
)}\frac{5}{2}\cdot
\frac{3}{2}h^4\bigl(|\xi|^2+\lambda_dh^2
\bigr)^{\sfrac
{1}{2}} \cdot\mathcal{F}\eta(\xi),
\end{eqnarray*}
respectively.
It can be shown by a straightforward calculation that $\xi^{\alpha
}f_j\in\mathcal{W}^{6+|\alpha|}(\R^d)$ for all $\alpha\in\{
0,\dots,d\}^d$.
\end{enumerate}
\end{Example}

\begin{Remark}\label{DiscAss}In the one-dimensional regression model
$\eqref{model}$, Birke \textit{et al.} \cite{birbishol2010} assume
that the kernel $K$ has exponentially decreasing tails in order to
obtain asymptotic confidence bands, which, in
combination with the other assumptions only allows for kernels that are
Fourier transforms of $C^{\infty}$-functions
with square integrable derivatives. Our Assumption \ref{assB} is already
satisfied if $K$ is the Fourier transform of a once
weakly differentiable function with square integrable weak derivative,
such that all indices of ill-posedness $\beta$
that satisfy $\beta>\frac{1}{2}$ are included if $d=1$. Moreover, the
assumptions regarding the bandwidths are less
restrictive compared to Birke \textit{et al.} \cite{birbishol2010}.
\end{Remark}

Our final assumptions refer to the smoothness of the function $f$ and
to the decay of the convolution $f*\psi$.

\begin{Assumption}\label{assC} We assume that
\begin{enumerate}[C2.]

\item[C1.] There exist constants $\gamma>2$, $m>\gamma+\frac{d}{2}$
such that $f\in\mathcal{W}^{m}(\R^d)$.
\item[C2.] There exists a constant $\nu>0$ such that
\begin{eqnarray*}
\int_{\R}\bigl\llvert (f*\psi) (z)\bigr\rrvert
^2\bigl(1+|z|^2\bigr)^{\nu} \,\mathrm{d}z<\infty.
\end{eqnarray*}
\end{enumerate}
\end{Assumption}

\section{Asymptotic confidence regions}\label{sec4}

In this section, we construct asymptotic confidence regions for the
function $f$ on the unit cube $[0,1]^d$. These results can easily be
generalized to arbitrary rectangles  $\bigtimes_{j=1}^d[a_j,b_j ]$ for
fixed constants $a_j<b_j$ $(j=1,\ldots,d)$ and the details are omitted
for the sake of brevity. We
investigate the limiting distribution of the supremum of the process $\{
\tilde Y_n(x) | x\in[0,1]^d\}$, where
%
\begin{eqnarray}
\label{DS} \tilde Y_n(x)&=&\frac{(2\pi)^dh^{\beta}\sqrt{h^d
n^{d}a_n^d}}{\sigma\|K\|_{L^2(\R^d)}} \bigl[\hat
f_n(x)-\E\hat f_n(x) \bigr]
\nonumber\\[-8pt]\\[-8pt]
&=&\frac{(2\pi)^dh^{\beta}}{\sigma\|K\|_{L^2(\R^d)}\sqrt{h^d
n^{d}a_n^d}}\sum_{\mathbf{k}\in G_n}K_n
\biggl( (x-x_{\mathbf
{k}} ) \frac{1}{h} \biggr)\varepsilon_{\mathbf{k}}\nonumber
\end{eqnarray}
and the kernel $K_n$ is defined in \eqref{kernKn}. Note that
\begin{eqnarray*}
\sup_{x\in[0,1]^d}\bigl|\tilde Y_n(x)\bigr|= \sup
_{x\in[0,h^{-1}]^d}\bigl| Y_n(x)\bigr|,
\end{eqnarray*}
where the process
%
\begin{eqnarray}
\label{yn} Y_n(x):= \frac{(2\pi)^dh^{\beta}}{\sigma\|K\|_{L^2(\R^d)}\sqrt {h^d n^{d}a_n^d}} \sum
_{\mathbf{k}\in G_n}K_n \biggl(x-x_{\mathbf{k}}
\frac{1}{h} \biggr)\varepsilon_{\mathbf{k}}
\end{eqnarray}
can be approximated by a stationary Gaussian field uniformly with
respect to
$[0,h^{-1}]^d$. Thus the desired limiting distribution corresponds to
the limiting distribution of the supremum of a stationary Gaussian
process over a system of increasing smooth sets with sufficient
similarity of their speed of increase, and is therefore of Gumbel-type.
The precise result is given in the following theorem.

\begin{Theorem}\label{GWS}
Assume that for some fixed constant $\delta\in(0,1]$, $\delta<d$ and
a constant $r>\frac{2d}{d-\delta}$ the $r$th moment of the errors
exists, that is, $\mathbb{E}|\varepsilon_{\mathbf{k}}|^r<\infty$.
If additionally Assumptions \textup{\ref{assA}} and \textup{\ref{assB}} are satisfied and $\frac{\log
(n)}{n^{\delta}a_n^{\delta}h^{d}}=\mathrm{o}(1)$, then we have
\begin{eqnarray*}
\lim_{n\rightarrow\infty}\Prob \Bigl(\sup_{x\in[0,1]^d} \bigl(\bigl|
\tilde Y_n(x)\bigr|-C_{n,3} \bigr)\cdot C_{n,3}<
\kappa \Bigr)=\mathrm{e}^{-2\mathrm{e}^{-\kappa}},
\end{eqnarray*}
where
\begin{eqnarray*}
C_1&=&\operatorname{\mathbf{det}} \biggl( \biggl[\frac{(2\pi)^{2d}}{\|K\|_2^2}\int
_{\R^d}\bigl|\Psi(v)\mathcal{F}\eta(v)\bigr|^2v_iv_j
\,\mathrm{d}v \biggr] , i,j=1,\ldots,d \biggr),
\\
C_{n,2}&=&\sqrt{\frac{C_1}{(2\pi)^{d+1}}}\frac{1}{h^d},
\\
C_{n,3}&=&\sqrt{2\ln(C_{n,2})}+
\frac{(d-1)\ln (2\ln
(C_{n,2}) )}{2\sqrt{2\ln(C_{n,2})}}.
\end{eqnarray*}
\end{Theorem}

The proof of this result is long and complicated and therefore deferred
to Sections~\ref{sec6} and \ref{sec7}.
In the following, we apply Theorem~\ref{GWS} to construct uniform
confidence regions for the function $f$ by choosing the bandwidth such
that the bias decays to zero sufficiently fast. More precisely, if the condition
\begin{eqnarray*}
\log(n) \sup_{x \in[0,1]^d}\bigl\llvert f(x)-\E\hat
f_n(x)\bigr\rrvert =\mathrm{o} \Bigl( \Bigl(h^{\beta}\sqrt
{h^dn^da_n^d}
\Bigr)^{-1} \Bigr)
\end{eqnarray*}
is satisfied, it follows directly that the random quantities $\sup_{x
\in[0,1]^d} |\tilde Y_n(x)|$ and
\begin{eqnarray*}
\frac{(2\pi)^dh^{\beta}\sqrt{h^dn^da_n^d}}{\|K\|_{L^2(\R^d)}\sigma
}\sup_{x\in[0,1]^d}\bigl\llvert f(x)-\hat
f_n(x)\bigr\rrvert
\end{eqnarray*}
have the same limiting behavior.

\begin{Corollary}\label{cor1}
Assume that the conditions of Theorem~\ref{GWS}, Assumption \textup{\ref{assC}} and the condition
\begin{eqnarray*}
\sqrt{h^dn^da_n^d}
\sqrt{\log(n)} \biggl(\frac{1}{n^3a_n^3h^2}+\frac
{a_n^{\nu}}{n}+a_n^{\nu+\sfrac{d}{2}}+h^{\gamma+\beta}
\biggr)=\mathrm{o}(1)\qquad  \mbox{for } n\rightarrow\infty
\end{eqnarray*}
are satisfied. Then we have for any $\kappa\in\R$
\begin{eqnarray*}
\lim_{n\rightarrow\infty}\Prob \bigl(\hat f_n(x)-
\Phi_{n,\kappa
}\leq f(x)\leq\hat f_n(x)+\Phi_{n,\kappa}
\mbox{ for all } x\in [0,1]^d \bigr)=\mathrm{e}^{-2\mathrm{e}^{-\kappa}},
\end{eqnarray*}
where the sequence $\Phi_{n,\kappa}$ is defined by
\begin{eqnarray*}
\Phi_{n, \kappa}=\frac{(\sfrac{\kappa}{C_{n,3}}+C_{n,3})\sigma\|K\|
_{L^2(\R^d)}}{(2\pi)^dh^{\beta}\sqrt{h^dn^da_n^d}}.
\end{eqnarray*}
\end{Corollary}

As a consequence of Corollary~\ref{cor1} an asymptotic uniform
confidence region for the function $f$ with confidence level $1-\alpha
$ is given by
%
\begin{eqnarray}
\label{conf} \bigl\{\bigl[\hat f_n(x)-\Phi_{n,-\ln(-0.5\ln(1-\alpha))},\hat
f_n(x)+\Phi _{n,-\ln(-0.5\ln(1-\alpha))}\bigr] | x\in [0,1]^d\bigr
\}.
\end{eqnarray}
The corresponding $(1-\alpha)$-band has a width of $2\Phi_{n,-\ln
(-0.5\ln(1-\alpha))}$. Here, the factor $\frac{1}{h^{\beta}}$ is
due to the ill-posedness of the inverse problem (see Assumption \ref{assB}). It
does not appear in corresponding results for the direct regression
case. On the other hand, the factor $a_n^{-\sfrac{d}{2}}$ arises from
the design on the growing system of sets $\{[-a_n^{-1},a_n^{-1}]^d |
n\in\N\}$. In the case of a regression on a fixed interval it does
not appear as well. The width of the asymptotic point-wise confidence
intervals in the multivariate indirect regression case as obtained in
Bissantz and Birke \cite{bisbir2009} is of order $\frac
{1}{h^{\beta}\sqrt {Nh^da_n^{d}}}$, where $N$ is the total number of observations. Their
point-wise confidence intervals are smaller than the uniform ones
obtained in Corollary~\ref{cor1}. The price for uniformity is an
additional factor of logarithmic order, which is typical for results of
this kind.

In applications the standard deviation is unknown but can be estimated
easily from the data, because this does not require the estimation of
the function $f$. In particular, \eqref{conf} remains an asymptotic
$(1- \alpha)$-confidence band, if $\sigma$ is replaced by an
estimator satisfying
$\hat\sigma-\sigma=\mathrm{o}_P(1/\log(n))$.

\section{Time dependent regression functions}\label{sec5}

In this section, we extend model $(\ref{model})$ to include a time
dependent regression function, that is
%
\begin{eqnarray}
\label{model2} Y_{j,\mathbf{k},n}= (T_{\psi}f_{t_j} )
(x_{\mathbf
{k}})+\varepsilon_{\mathbf{k}}, \qquad \mathbf{k}\in G_{\mathbf{n}},
j=-m,\ldots,m,
\end{eqnarray}
where $x_{\mathbf{k}}=\frac{\mathbf{k}}{na_n}$ and $t_j=\frac{j}{mb_m}$,
$m=m(n)$, such that $m(n)\rightarrow\infty$ and $b_{m(n)}\searrow0$
as $n\rightarrow\infty$.

We assume that $\psi$ does not depend on the time and the operator
$T_{\psi}$ is defined by
\begin{eqnarray*}
(T_{\psi}f_{t} )=\int_{\mathbb{R}^d}f_t(y)
\psi(\cdot-y) \,\mathrm{d}y.
\end{eqnarray*}
This assumption is reasonable in the context of imaging where the
function $\psi$ corresponds to the point spread function (Bertero \textit{et al.} \cite{berbocdes2009}).
If it is not satisfied, that is, the convolution operator effects all
coordinates, the problem can be modeled as in Section~\ref{sec3}.

For a precise statement of the results, we will add an index to the
Fourier operator $\mathcal{F}$ which gives the dimension of the space
under consideration. We will write $\mathcal{F}_{d+1}$ if the Fourier
transform is taken over the whole space $\R^{d+1}$ and $\mathcal
{F}_d$ to denote
Fourier transformation with respect to the spatial dimensions.
By the same considerations as given in Section~\ref{sec3}, we obtain
an estimator $\check f$ for the function $f_t$
\begin{eqnarray*}
\check f_n(x;t)&=&\frac{1}{(2\pi)^{\vfrac{d+1}{2}}}\int_{\R^{d+1}}
\frac{\mathcal{F}_{d+1}(\widehat{f*\psi})(\xi,\tau)}{(2\pi
)^{\sfrac{d}{2}}\mathcal{F}_{d}\psi(\xi)} \mathcal{F}_{d}\check \eta(\xi h,\tau
h_t)\exp\bigl(\mathrm{i}t\tau+\mathrm{i}x^T\xi\bigr) \,\mathrm{d}(\xi,\tau)
\\
&=&\frac{1}{(2\pi)^{d+\sfrac{1}{2}}n^dma_n^db_m}\sum_{(\mathbf
{k},j)\in G^{d+1}_{(\mathbf{n},m)}}Y_{\mathbf{k},j}
\check K_n \biggl( \frac{x-x_{\mathbf{k}}}{h},\frac{t-t_j}{h_t} \biggr),
\end{eqnarray*}
where $G^{d+1}_{(n,m)}$ denotes the grid $\{-n,\ldots,n\}^d\times\{
-m,\ldots,m\}$ and the kernel $\check K_n$ is given by
%
\begin{eqnarray}
\label{KernZeitUngefaltet} \check K_n(x;t)=\frac{1}{(2\pi)^{\vfrac{d+1}{2}}}\int
_{\R^{d+1}} \frac{\exp (\mathrm{i}\tau t+\mathrm{i}\xi^Tx )}{\mathcal{F}_{d}\psi(\sfrac
{\xi}{h})}\mathcal{F}_{d+1}\check\eta(
\xi,\tau) \,\mathrm{d}(\xi,\tau).
\end{eqnarray}
Here the function $\check\eta\dvtx \R^{d+1}\rightarrow\R$ satisfies
condition \ref{assA} and $h_t=h_t(n)$ is an additional sequence of bandwidths
referring to the time domain.
For the asymptotic analysis, we require a modified version of
Assumption \ref{assB}.

{\renewcommand{\theAssumption}{$\bolds{\check B}$}
\begin{Assumption}\label{assBB}
Let Assumptions \textup{B1} (with corresponding kernel $\check K$) and \textup{B2} hold
and additionally assume that
\begin{enumerate}[$\check \mathrm{B}3$.]
\item[$\check \mathrm{B}3$.] $\log(n+m(n))\cdot h^{\mu_M}(a_n^{-\sfrac
{d}{2}}h^{-\sfrac{d}{2}}b_{m(n)}^{\sfrac{1}{2}}m(n)^{\sfrac
{1}{2}})=\mathrm{o}(1)$ and for
$p=1,\ldots,M-1$
\[
h^{\mu_{p}}\bigl(\log(n+m)\bigr)^2=\mathrm{o}(1).
\]
\end{enumerate}
\end{Assumption}}

\begin{Theorem}\label{GWSZeit}
{Define}
\begin{eqnarray*}
\check Y_n(x;t):=\frac{(2\pi)^{d+1}h^{\beta}\sqrt{h^dh_t
n^{d}mb_ma_n^d}}{\sigma\|\check K\|_{L^2(\R^{d+1})}\|} \bigl[\check
f_n(x;t)-\E\check f_n(x;t) \bigr]
\end{eqnarray*}
and let the moment condition of Theorem~\ref{GWS} and Assumptions \textup{\ref{assA}}
and \textup{\hyperref[assBB]{$\check \mathrm{B}$}} be satisfied.
We further assume that the bandwidths $h_t$ and $h$, and the sequences
$(a_n)_{n\in\N}$ and $(b_{m(n)})_{n\in\N}$ satisfy
\begin{eqnarray*}
&\displaystyle\log(n+m) \biggl(\sqrt{\frac{na_n}{mb_m}}\frac{1}{\sqrt{n^{\delta
}h_ta_n^{\delta}h^d}}+ \biggl(
\frac{mb_m}{na_n} \biggr)^{\sfrac
{d}{2}}\frac{1}{\sqrt{m^{\delta}h_th^{d}}} \biggr) =\mathrm{o}(1) \qquad \mbox{for } n\rightarrow\infty,&
\\
&\displaystyle h_t+h\leq L_1\cdot h^{d(1-L_2)}h_t^{(1-L_2)}&
\end{eqnarray*}
for some constants $L_1<\infty$ and $L_2 \in(0,1)$.
Then we have for each $\kappa\in\R$,
\begin{eqnarray*}
\lim_{n\rightarrow\infty}\Prob \Bigl(\sup_{x\in[0,1]^d} \bigl(\bigl|
\check Y_n(x;t)\bigr|-D_{n,3} \bigr)\cdot D_{n,3}<
\kappa \Bigr) =\mathrm{e}^{-2\mathrm{e}^{-\kappa}},
\end{eqnarray*}
where
\begin{eqnarray*}
D_1&=&\operatorname{\mathbf{det}} \biggl( \biggl[\frac{(2\pi)^{2(d+1)}}{\|\check K\|
^2_{L^2(\R^{d+1})}}\int
_{\R^{d+1}}\bigl|\Psi(v_1,\ldots,v_d)
\mathcal{F}_{d+1} \check\eta(v)\bigr|^2v_iv_j
\,\mathrm{d}v \biggr] , i,j=1,\ldots,d+1 \biggr),
\\
D_{n,2}&=&\sqrt{\frac{D_1}{(2\pi)^{d+2}}}\frac{1}{h^dh_t}
\quad \mbox{and}
\\
D_{n,3}&=&\sqrt{2\ln(D_{n,2})}+
\frac{(d-1)\ln (2\ln
(D_{n,2}) )}{2\sqrt{2\ln(D_{n,2})}}.
\end{eqnarray*}
\end{Theorem}

\begin{Corollary} \label{cor2}
If the assumptions of Theorem~\ref{GWSZeit} are satisfied, the limit
kernel $\check K$ is defined by
%
\begin{eqnarray}
\label{kcheck} \check K(x,t)=\frac{1}{(2\pi)^{\vfrac{d+1}{2}}}\int_{\R^{d+1}}
\Psi (\xi)\mathcal{F}_{d+1}\check\eta(\xi,\tau)\exp\bigl(\mathrm{i}
\xi^T x+\mathrm{i}\tau t\bigr)\, \mathrm{d}(\xi,\tau)
\end{eqnarray}
and the function $f_{(\cdot)}(\cdot)\dvtx \R^{d+1}\rightarrow\R^1$
satisfies Assumption \textup{\ref{assC}}, then it follows that
\begin{eqnarray*}
\lim_{n\rightarrow\infty}\Prob \bigl(\check f_n(x;t)-\check
\Phi _{n,\kappa}\leq f(x;t)\leq\check f_n(x;t)+\check
\Phi_{n,\kappa} \mbox{ for all } (x,t)\in[0,1]^{d+1}
\bigr)=\mathrm{e}^{-2\mathrm{e}^{-\kappa}},
\end{eqnarray*}
where the constant $\check\Phi_{n,\kappa}$ is defined by
\begin{eqnarray*}
\check\Phi_{n,\kappa}=\frac{(\sfrac{\kappa
}{D_{n,3}}+D_{n,3})\sigma\|\check K\|_{L^2(\R^{d+1})}}{h^{\beta
}\sqrt{h^dn^da_n^dmb_mh_t}(2\pi)^{d+1}}.
\end{eqnarray*}
\end{Corollary}

Asymptotic confidence bands for the function $f_t(x)$ at level
$1-\alpha$ are hence given by
\begin{eqnarray*}
\bigl\{\bigl[\check f_n(x;t)-\check\Phi_{n,-\ln(-0.5\ln(1-\alpha))},\check
f_n(x;t)+\check\Phi_{n,-\ln(-0.5\ln(1-\alpha))}\bigr] | (x,t)\in
[0,1]^{d+1}\bigr\}.
\end{eqnarray*}

\section{Finite sample properties} \label{sec:sim}
In this section, we investigate the finite sample properties of the
proposed asymptotic confidence bands by means of a small simulation
study and
illustrate the procedure in an example analyzing the luminosity profile
of the elliptical galaxy NGC5017.

%
\subsection{Simulation study}\label{sec:setup}
\begin{figure}

\includegraphics{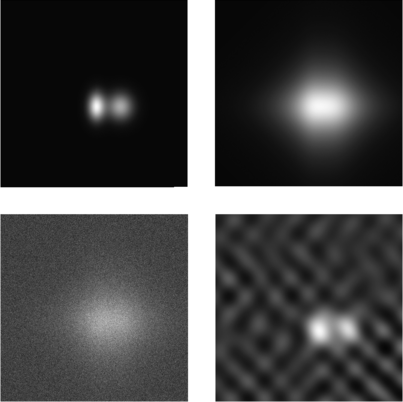}

\caption{Upper panel:
bimodal regression function $f_2$ (left) and convolution $f_2*\psi$
(right) defined in \protect\eqref{f2}.
Lower panel: simulated data from $f_2$ observations according to the
model (\protect\ref{mod1}) (left) and the corresponding reconstruction
of $f_2$ (right). (See Section \protect\ref{sec:setup} for details of the
choice of the (slightly undermoothing) bandwidth.)}\label{fig:vier}\vspace*{-3pt}
\end{figure}

All results are based on $5000$ simulation runs.
We simulate data from the bivariate regression model \eqref{mod1}
where the errors $\varepsilon_{(k_1,k_2)}$ are independent and
$\mathcal{N}(0,\sigma^2)$-distributed, where $\sigma=0.1,
(k_1,k_2)\in\{-n,\ldots,n\}^2, n\in\{150,250,350,500,650\}$ and our
dataset provides $n=150$, which corresponds to a grid of $301\times
301$ data points. For the unknown regression function we consider both,
a unimodal function
\begin{eqnarray*}
f_1(x,y)=4\exp \bigl(-(3.5x-1.5)^2-(3.5y-1.5)^2
\bigr)
\end{eqnarray*}
and a bimodal function
%
\begin{equation}
\label{f2} \hspace*{-5pt}f_2(x,y)=4\exp\bigl(-\bigl(3(x-0.1)
\bigr)^2-3.5(y-0.5)^2\bigr)+3\exp\bigl(-\bigl(2(x-1)
\bigr)^2-3.5(y-0.5)^2\bigr).
\end{equation}
As convolution function $\psi$ we consider $\psi(x,y)=\frac
{1}{4}\exp (-(|x|+|y|) )$ and the values for the sequence
$(a_n)_{n\in\N}$ are chosen such that most of the signal considerably
different from $0$ is included in the observations, that is
$(a_{150},a_{250},a_{350},a_{500},a_{650})=(0.29,0.28,0.27,0.26,0.25)$.
A difference-based variance estimator is used to estimate $\sigma$.
Figure~\ref{fig:vier} shows exemplary one simulated dataset and the
reconstruction of the bimodal regression function $f_2$ from this
dataset in comparison to the function $f_2$ itself and the convolution
$f_2*\psi$.
%
%
\begin{table}[b]\vspace*{-3pt}
\tablewidth=\textwidth
\tabcolsep=0pt
\caption{Simulated coverage probabilities and
mean half-lengths of the corresponding confidence bands for the
bivariate Gaussian function $f_1$ and the bimodal function $f_2$}
\label{table}
\begin{tabular*}{\textwidth}{@{\extracolsep{4in minus 4in}}llll@{\hspace*{-16pt}}lll@{\hspace*{-16pt}}lll@{\hspace*{-16pt}}ll@{}} \hline
                                                                                 &        &       & \multicolumn{3}{l}{$90\%$ nominal coverage}           & \multicolumn{3}{l}{$95\%$ nominal coverage}           & \multicolumn{3}{l}{$99\%$ nominal coverage}           \\
 [-5pt]
                                                                                 &        &       & \multicolumn{3}{l}{\hrulefill}                        & \multicolumn{3}{l}{\hrulefill}                        & \multicolumn{3}{l}{\hrulefill}                        \\
$n$                                                                              & $a_n$  & $f$   & 100-Cov. $(\%)$                      &           & Length  & 100-Cov. $(\%)$                      &           & Length  & 100-Cov. $(\%)$            &                     & Length  \\\hline
$150$                                                                            & $0.29$ & $f_1$ & $11.5$                           &           & $0.235$ & $3.6$                            &           & $0.240$ & $ 0.1$                 &                     & $0.285$ \\
                                                                                 &        & $f_2$ & $\phantom{1}{9.0}$               &           & $1.219$ & $4.4$                            &           & $0.324$ & $0.8 $                 &                     & $0.368$ \\[3pt]
$250$                                                                            & $0.28$ & $f_1$ & $10.6$                           &           & $0.159$ & $3.9$                            &           & $0.169$ & $0.2 $                 &                     & $0.192$ \\
                                                                                 &        & $f_2$ & $\phantom{1}{9.7}$               &           & $0.265$ & $4.8$                            &           & $0.290$ & $0.6 $                 &                     & $0.319$ \\[3pt]
$350$                                                                            & $0.27$ & $f_1$ & $10.0$                           &           & $0.109$ & $5.0$                            &           & $0.146$ & $0.5 $                 &                     & $0.166$ \\
                                                                                 &        & $f_2$ & $10.4$                           &           & $0.240$ & $5.9$                            &           & $0.254$ & $ 0.6$                 &                     & $0.287$ \\
[3pt]
$500$                                                                            & $0.26$ & $f_1$ & $\phantom{1}{8.6}$               &           & $0.108$ & $4.3$                            &           & $0.115$ & $0.6 $                 &                     & $0.130$ \\
                                                                                 &        & $f_2$ & $\phantom{1}{9.1}$               &           & $0.191$ & $5.4$                            &           & $0.203$ & $0.8 $                 &                     & $0.229$ \\[3pt]
$650$                                                                            & $0.25$ & $f_1$ & $10.3$                           &           & $0.412$ & $5.3$                            &           & $0.437$ & $0.6 $                 &                     & $0.495$ \\
                                                                                 &        & $f_2$ & $\phantom{1}{9.7}$               &           & $0.648$ & $4.8$                            &           & $0.687$ & $0.6 $                 &                     & $0.775$ \\
\hline
\end{tabular*}
\end{table}

\begin{figure}

\includegraphics{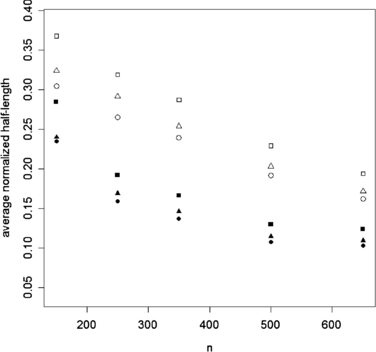}

\caption{Confidence band average half-lengths for both regression
functions, $f_1\  (\alpha=0.01\dvt \blacksquare, \alpha
=0.05$:  $\blacktriangle, \alpha=0.1\dvt \mbox{\CIB})$ and $f_2\ (\alpha
=0.01\dvt \square, \alpha=0.05\dvt \vartriangle, \alpha=0.1\dvt \mbox{\CIW})$.}
\label{fig:coverage}
\end{figure}

\begin{figure}[b]

\includegraphics{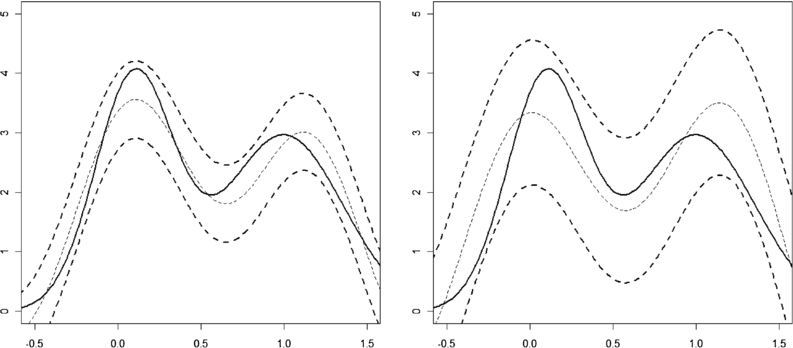}

\caption{Cross-section of bivariate, bimodal regression function $f_2$
(solid line), cross-section of reconstructions of $f_2$ (thin dashed
lines) and corresponding 90\%-confidence bands (thick dashed lines) for
$n=650$ (left) and $n=150$ (right).}
\label{fig:Baender}
\end{figure}

For computational feasibility, we determine at first for each scenario
a bandwidth by a small preliminary simulation study.
For this purpose, we applied the $L_{\infty}$-motivated bandwidth
selection method introduced in Bissantz \textit{et al.} \cite
{bisdumholmun2007}
and the estimated bandwidth is used in all $5000$ runs.
Table~\ref{table} shows the simulated coverage probabilities and the
average half-widths of the bands normalized with respect to the maximum
of the respective function. Figure~\ref{fig:coverage} illustrates the
decrease of the normalized average half-widths of the confidence bands
plotted against the sample sizes for both, the unimodal and the bimodal
function. We conclude that for larger sample sizes (note that $n=150$
corresponds to $301\times301$ observations) and relatively small
variances the simulated coverage probabilities are close to the nominal
values. For $n=150$ the bands are rather wide, especially for the
regression function $f_2$, that requires a smaller bandwidth for the
estimation than the function $f_1$, which results in considerably wider
bands for the function $f_2$ than for $f_1$. For increasing sample
sizes, the widths of the bands decrease significantly. For
illustrational purposes,
Figure~\ref{fig:Baender} shows a cross-section of the bivariate
function $f_2$, estimators and a set of $90\%$ confidence bands for
$n=150$ and $n=650$ respectively,
where $y=0.5$ and the corresponding confidence bands have been obtained
in the bivariate setting. This figure clearly demonstrates the increase
of precision of the bands for increasing sample size. Note that the
sample size $301\times301$ in this example
is rather small compared to the sample sizes which are usually
available for astronomical images.
Moreover, in these applications the signal-to-noise ratio is often much
smaller and the point spread function is usually more sharply
peaked than the one used in the simulations.
%
%
\subsection{The luminosity profile of the galaxy NGC5017}

\begin{figure}

\includegraphics{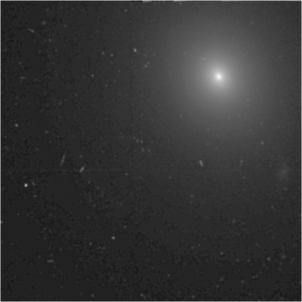}

\caption{Hubble Space Telescope/Wide Field Planetary Camera 2
[HST/WFPC2] R-band image of the
elliptical galaxy NGC5017. Image source: NASA/STScI, dataset U3CM1A01R.}\label{fig:Galaxie}
\end{figure}

In this section, we use the methodology derived in this paper to
analyze the shape of luminosity profiles of elliptical galaxies.
Figure~\ref{fig:Galaxie} shows a $301\times301$ pixel section of an
HST/WFPC2 R-band image of the elliptical galaxy NGC5017.
%

It is well known that images, taken with telescopes, are usually at
least slightly blurry which can be modelled as convolution of the sharp
image with a so-called point spread function (PSF), $\psi$, of the
optical instrument. The detector used is a digital imaging device (CCD,
charge-coupled device). We use a dataset of the size $301\times301$, where
each data point corresponds to a pixel on an equally spaced grid.
Hence, the two dimensional model \eqref{mod1} is suitable to describe
the data.

\begin{figure}

\includegraphics{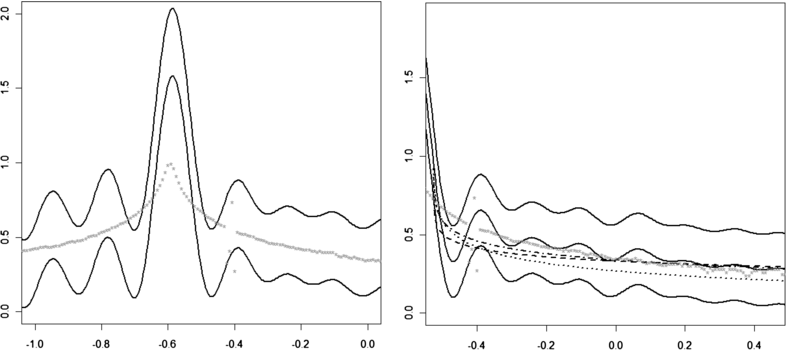}

\caption{Upper and lower $90\%$-confidence bands for the luminosity
profile of the galaxy (solid lines) and corresponding data points
($\ast$).
The left plot shows the central region
and the right plot the region in which the profile is fitted with S\'
ersic-curves for $\kappa=4$ (dotted line), $\kappa=5.11$
(dashed-dotted line) and $\kappa=7$ (dashed line).}
\label{fig:cross}
\end{figure}

In the analysis of elliptical galaxies the luminosity profile, that is,
the decrease of the brightness with increasing distance from the
galactic centre, is of particular interest. We use the confidence bands
based on the nonparametric estimator \eqref{est} to narrow down the
region for the parameter $\kappa$ in the S\'ersic (1968) model for the
luminosity profile of NGC5017. This model is defined as
%
\begin{eqnarray}
\label{sersic} f^{\kappa}(r)=I_0\cdot\exp
\bigl(-b_{\kappa}(r/r_e)^{\sfrac
{1}{\kappa}} \bigr),
\end{eqnarray}
where $r$ is the distance from the centre of the galaxy, $I_0$ is the
brightness in the centre $(r=0)$, $r_e$ is the scale radius (i.e., the
half-light radius) and $b_{\kappa}$ is a normalization constant that
is uniquely determined by the choice of $\kappa$ which is the shape
parameter, controlling the curvature of the profile. Note that $I_0$
and $r_e$ are model independent quantities that can be found in the
literature. For $\kappa=4$, model \eqref{sersic} coincides with the
famous de Vaucouleur-profile (1959).
For details see, for example, Trujillo \textit{et al.} \cite{truj} who
already classified the
galaxy under consideration as S\'ersic-type galaxy with $\kappa=5.11$.
To analyse the data, we first fit the PSF, given by
\begin{eqnarray*}
\psi(x,y)=\frac{\lambda}{2}\exp\bigl(-\lambda\bigl(|x|+|y|\bigr)\bigr)
\end{eqnarray*}
to the (probably stellar) point source on the left middle of the image
since the PSF is not completely known in advance. The bandwidth was
chosen according to the procedure described in Section~\ref{sec:setup}. We compute the estimator and $90\%$-confidence bands for
the unknown luminosity profile which, in combination, suggest that the
 actual brightness gradient is clearly steeper in the central region
than the (convolved), unprocessed data tell us, see Figure~\ref{fig:cross}, left panel, for a cross-section. The oscillations are
artifacts due to our use of a Fourier estimator whose severity is
partly caused by the use of one constant bandwidth for the
reconstrucion of the whole image.
In order to reconstruct the central region of the galaxy with
sufficient precision, the bandwidth is chosen too small for the
shallower regions of the image. For the analysis, we restrict ourselves
to the square $[-0.55,0.45]^2$ and check for which parameters $\kappa
\in[4,8]$ the resulting profile $f^{\kappa}$ according to formula
\eqref{sersic} is completely contained in the
90\%-confidence band. We find that this is satisfied for the parameters
$\kappa\in[4.81,5.93]$, which suggests that the profile for $\kappa
=4$, corresponding to the de Vaucouleur law is not appropriate to
describe the data well. On the other hand, the data do not provide
evidence against $\kappa=5.11$, as proposed by Trujillo \textit{et al.}
\cite{truj}.

%
%

\section{Proofs of Theorem 
\texorpdfstring{\protect\ref{GWS}}{1} and Corollary \texorpdfstring{\protect\ref{cor1}}{1}}
\label{sec6}
\setcounter{equation}{0}

\subsection{Notation, preliminaries and remarks}\label{PreProof}
First, we introduce some notation which is used extensively in the
following proofs. Define for $a=(a_1,\dots,a_d),b=(b_1,\dots,b_d)\in
\R^d$ the $d$-dimensional cube $[a,b]:=\bigtimes_{j=1}^d[a_j,b_j]$.
Let $\mathbf{k}=(k_1,\ldots,k_d)\in\mathbb{Z}^d$, $\alpha= (\alpha
_1,\dots,\alpha_d) \in\{0,1\}^d$ be multi-indices, $\mathbf
{0}:=(0,\ldots,0)^T\in\R^d$ and $\mathbf{1}:=(1,\ldots,1)^T\in\R
^d$ and
define $G_{\mathbf{k}}:=\mathbb{Z}^d\cap[-\mathbf{k}, \mathbf
{k}]$. For $j \in\{ 1,\dots,d \}$ we denote by
$G_{\mathbf{k}}^j$ the canonical projection of $G_{\mathbf{k}}$ onto
$\mathbb{Z}^j$, i.e., $G_{\mathbf{k}}^j$ is a $j$-dimensional grid of
integers with possibly different length in each direction. For $j\in\N
$ let $G_{\mathbf{k}}^{j,+}:= G_{\mathbf{k}}^j\cap\N^j$ denote the
part of the grid $G_{\mathbf{k}}^j$ whose vectors contain only
positive components and write $G_{\mathbf{k}}^+$ for $G_{\mathbf
{k}}^{d,+}$. We further introduce the bijective map\vspace*{-1pt}
\begin{eqnarray*}
E_d\dvtx \lleft\{ %
\begin{array} {l} \{0,1\}^d
 \rightarrow \mathcal{P}\bigl(\{1,\ldots,d\}\bigr) ,
\\\noalign{\vspace*{-1pt}}
(\alpha_1,\ldots,\alpha_d)  \mapsto v=
\{v_1,\ldots,v_{|\alpha
|}\}; \qquad  \alpha_{v_j}=1, j=1,
\ldots,|\alpha| =\displaystyle \sum_{i=1}^d
\alpha_i, \end{array} %
\rright.
\end{eqnarray*}
that maps each $\alpha$ to the set $v\subset\{1,\ldots,d\}$ that
contains the positions of its ones. For $\alpha\in\{0,1\}^d$ and $\{
v_1,\ldots,v_{|\alpha|}\}=E_d(\alpha)$ let $(x)_{\alpha
}:=(x_{v_1},x_{v_{2}},\ldots,x_{v_{|\alpha|}})$ denote the projection
of $x \in\mathbb{R}^d$ onto the space spanned by the coordinate axes
given by the positions of ones of the multi-index $\alpha$. For
$a,b\in\R^d$ let $(a)_{\alpha}\dvtx (b)_{\mathbf{1}-\alpha
}=(a\dvtx b)_{\alpha}:=(a_1^{\alpha_1}\cdot b_1^{1-\alpha_1},\ldots
,a_d^{\alpha_d}\cdot b_d^{1-\alpha_d})$
denote the vector of the components of $a$ and $b$ specified by the
index $\alpha$. The following example illustrates these notations.

\begin{Example}For $d=2$ we have $\{0,1\}^2=\{(1,1),(1,0),(0,1),(0,0)\}
$ and the mapping $E_2$ is defined by\vspace*{-2pt}
\[
E_2 \bigl((1,1) \bigr)=\{1,2\},\qquad  E_2 \bigl((1,0) \bigr)=
\{1\},\qquad  E_2 \bigl((0,1) \bigr)=\{2\} \quad \mbox{and}\quad  E_2
\bigl((0,0) \bigr)=\emptyset.
\]
For any $x = (x_1,x_2) \in\R^2$ we have\vspace*{-2pt}
\[
(x)_{(1,1)}=x,\qquad  (x)_{(1,0)}=x_1 \quad \mbox{and}\quad
(x)_{(0,1)}=x_2.
\]
For $a=(a_1,a_2), b=(b_1,b_2)\in\R^2$ we have\vspace*{-2pt}
\begin{eqnarray*}
(a\dvtx b)_{(1,1)}&=&(a_1,a_2)=a, \qquad (a\dvtx b)_{(1,0)}=(a_1,b_2),\\[-1pt]
(a\dvtx b)_{(0,1)}&=&(b_1,a_2),\qquad  (a\dvtx b)_{(0,0)}=(b_1,b_2)=b.
\end{eqnarray*}
\end{Example}

For the approximation of the integrals by Riemann sums we define for
multi-indices $\tilde\alpha,\alpha\in\{0,1\}^d\backslash\{\mathbf
{0}\}$\vspace*{-2pt}
%
\begin{eqnarray}
\label{Delta} \Delta_{\alpha}(f;a,b) &:=&\sum_{\tilde\alpha\in\{0,1\}^d, \tilde\alpha\leq\alpha}
(-1)^{|\tilde\alpha|}f\bigl((a\dvtx b)_{\tilde\alpha}\bigr) \nonumber\\[-8pt]\\[-9pt]
&=&\sum
_{\tilde\alpha\in\{0,1\}^d, \tilde\alpha\leq\alpha} (-1)^{d-|\tilde\alpha|}f \bigl((a)_{\mathbf{1}-\tilde\alpha
}\dvtx (b)_{\tilde\alpha}
\bigr),\nonumber
\end{eqnarray}
where the symbol $\tilde\alpha\leq\alpha$ means $\tilde\alpha
_j\leq\alpha_j$ for $j=1,\ldots,d$.
Note that for $\alpha=\mathbf{1}\in\R^d$ we obtain the special case
of the $d$-fold alternating sum, that is,
\begin{eqnarray*}
\Delta(f;a,b) := \Delta_\mathbf{{1}} (f;a,b)= \sum
_{\alpha\in\{0,1\}
^d}(-1)^{|\alpha|}f \bigl((a\dvtx b)_{\alpha} \bigr)
=\sum_{\alpha\in\{0,1\}^d}(-1)^{d-|\alpha|}f
\bigl((a\dvtx b)_{\mathbf
{1}-\alpha}\bigr).
\end{eqnarray*}

Note that $\Delta_{\alpha}(f;a,b)$ can be regarded as the increment
of the function
$f_\alpha((x)_\alpha):= f ((x\dvtx b)_{\alpha} )$
over the interval $[(a)_{\alpha},(b)_{\alpha}]$
which also gives rise to the alternative notation
%
\begin{eqnarray}
\label{Zuwachs} \Delta_{\alpha} (f;a,b )=\Delta \bigl(f_{\alpha
},(a)_{\alpha},(b)_{\alpha}
\bigr). 
\end{eqnarray}

\subsection{Proof of Theorem \texorpdfstring{\protect\ref{GWS}}{1}}

To prove the assertion of Theorem~\ref{GWS} we decompose the index set
$G_{\mathbf{n}}= \{-n,\ldots, n\}^d$ of the sum in $(\ref{DS})$ into
$2^d+1$ parts: the respective intersections with the $2^{d}$ orthants
of the origin in $\R^d$ and the marginal intersections with the
coordinate axes. Our first auxiliary result shows that
the contribution of the term representing the marginals is negligible
(here and throughout the paper we use the convention $0^0=1$).

\begin{Lemma} \label{Rand}
\begin{eqnarray*}
\sup_{x\in[0,h^{-1}]^d} \biggl|\frac{h^{\beta}}{\sqrt {h^dn^da_n^d}}\sum
_{\alpha\in\{0,1\}^{d}\backslash\{\mathbf{1}\}} \sum_{(\mathbf{k}:\mathbf{0})_{\alpha},\mathbf{k}\in
G^{+}_n}K_n
\biggl(x-\frac{1}{h}x_{\mathbf{k}} \biggr)\varepsilon _{\mathbf{k}}
\biggr|=\mathrm{o}_P \biggl(\frac{1}{\log(n)} \biggr).
\end{eqnarray*}
\end{Lemma}

We obtain from its proof in Section~\ref{sec7} that Lemma~\ref{Rand}
holds under the weaker condition
$\frac{\log(n)}{\sqrt{na_nh}}=\mathrm{o}(1)$, which follows from the
assumptions of Theorem~\ref{GWS}.
Next we consider the ``positive'' orthant $G_n^+$ and
show in three steps that
%
\begin{equation}
\label{approximation} \sup_{x\in[0,h^{-1}]^d} \bigl| Y_{n}^{(+)}(x)
- Y^{(+)}(x) \bigr| = \mathrm{o}_p(1),
\end{equation}
where the processes $ Y_{n}^{(+)}$ and $Y^{(+)}$ are defined by
%
\begin{eqnarray}
\label{Yn+} Y_{n}^{(+)}(x) &:= & \frac{(2\pi)^dh^{\beta}}{\sigma\|K\|
_{L^2}\sqrt{h^dn^da_n^d}} \sum
_{\mathbf{k}\in G_n^+} K_n\biggl(x-
\frac{1}{h}x_{\mathbf
{k}}\biggr)\varepsilon_{\mathbf{k}} ,
\\
\label{Y+} Y^{(+)}(x)&:=& \frac{(2\pi)^d}{\|K\|_{L^2}}\int
_{\R^d_+}K(x-u) \,\mathrm{d} B(u),
\end{eqnarray}
respectively,
$ B$ is a standard Brownian sheet on $\R^d$ (see the proof of
Lemma~\ref{Rio} for details) and
$K$ denotes the kernel defined in Assumption \ref{assB}.
The final result is then derived using Theorem~14.1 in Piterbarg \cite{piterbarg1996}.
To be precise note that it can easily be shown that
\begin{eqnarray*}
\lim_{n\rightarrow\infty}n^da_n^dh^dh^{2\beta}
\cdot\Var \bigl(\hat f_n(x) \bigr) =\frac{\sigma^2}{(2\pi)^{2d}}\int
_{\R^d}\biggl|K\biggl(\frac{x}{h}-u\biggr)\biggr|^2 \,\mathrm{d}u=
\frac{\sigma^2\|K\|^2_{L_2}}{(2\pi)^{2d}}
\end{eqnarray*}
(in particular the limit is independent of the variable $x$, which is
typical for kernel estimates in homoscedastic regression models with
equidistant design). We further obtain for the function
$ r(t)=(2\pi)^{2d}\|K\|^{-2}_{L^2}\int_{\R^d}K(v+t)K(v)\,\mathrm{d}v $ that
\begin{eqnarray*}
\|r\|_{L^1}=\frac{(2\pi)^{2d}}{\|K\|^2_{L^2}}\int_{\R^d}\biggl
\llvert \int_{\R^d}K(v+t)K(v) \,\mathrm{d}v\biggr\rrvert \,\mathrm{d}t\leq
\frac{(2\pi)^{2d}\|K\|
_{L^1}^2}{\|K\|^2_{L^2}}<\infty.
\end{eqnarray*}
Therefore the conditions of Theorem~14.1 in Piterbarg
\cite{piterbarg1996} are
satisfied and the assertion of Theorem~\ref{GWS} follows.

The remaining proof of the uniform approximation \eqref{approximation}
will be accomplished showing
the following auxiliary results. For this purpose we introduce the process
\begin{eqnarray*}
Y_{n,1}^{(+)}(x):=\frac{(2\pi)^dh^{\beta}}{\sqrt{n^da_n^dh^d}\|K\|
_{L^2(\R^d)}} \sum
_{\alpha\in\{0,1\}^d}(-1)^{|\alpha|}\sum_{\mathbf{j}\in
G_n^{|\alpha|,+}}
\Delta_{\alpha} (K_n\circ\tau_x,I_{\mathbf
{j}}
) B\bigl(\mathbf{j}\dvtx (\mathbf{n})_{\mathbf{1}-\alpha}\bigr),
\end{eqnarray*}
where the function $\tau_x$ is defined
by $\tau_x(u):=x-\frac{u+1}{na_nh}$,
%
\begin{eqnarray}
\label{Interval} I_{\mathbf{j}}:= \bigl[(\mathbf{j}-\mathbf{1})\dvtx (\mathbf
{n})_{1-\alpha}, \mathbf{j}\dvtx (\mathbf{n})_{1-\alpha} \bigr]\subset
\R_+^d
\end{eqnarray}
and we use the notation \eqref{Zuwachs}.

\begin{Lemma}\label{Rio}
There exists a Brownian sheet $ B$ on $\R^d$ such that
\begin{eqnarray*}
\sup_{x\in[0,h^{-1} ]^d}\bigl|Y_{n}^{(+)}(x)-Y_{n,1}^{(+)}(x)\bigr|=\mathrm{o}
\biggl(\frac{1}{\sqrt{\log(n)}} \biggr) \qquad \mbox{a.s.}
\end{eqnarray*}
\end{Lemma}

We obtain from the proof in Section~\ref{sec7.1} that Lemma~\ref{Rio}
holds under the condition $\frac{\log(n)}{n^{\sfrac{\delta
}{2}}a_n^{\sfrac{\delta}{2}}h^{\sfrac{d}{2}}}=\mathrm{o}(1)$, which
follows from the assumptions of Theorem~\ref{GWS}.
The next step consists of the replacement of the kernel $K_n$ in the
process $Y_{n,1}$ by its limit.

\begin{Lemma}\label{K_nK}
\begin{eqnarray*}
\sup_{x\in[0,h^{-1}]} \bigl|Y_{n,1}(x)-Y_{n,2}(x)
\bigr|=\mathrm{o}_P \biggl(\frac{1}{\log(n)} \biggr),
\end{eqnarray*}
where the process $Y_{n,2}$ is given by
\[
Y_{n,2}(x):=\frac{(2\pi)^d}{\sqrt{n^da_n^dh^d}\|K\|_{L^2}} \sum_{\alpha,\gamma\in\{0,1\}^d}(-1)^{|\alpha|}
\sum_{\mathbf
{j}\in
G_n^{|\alpha|,+}} \Delta_{\alpha} \bigl(K\circ
\tau_x,(-1)^{\gamma}I_{\mathbf
{j}} \bigr) B
\bigl((-1)^{\gamma}\mathbf{j}\dvtx (\mathbf{n})_{\mathbf{1}-\alpha}\bigr).
\]
\end{Lemma}

As described in Section~\ref{PreProof} for fixed $\alpha\in\{0,1\}
^{d}, j\in G_n^{|\alpha|,+}$ the quantity
$\Delta_{\alpha} (K\circ\tau_x; I_{\mathbf{j}} )$ can
be regarded as the increment of the function $(K_n\circ\tau
_x)_{\alpha}((u)_{\alpha})= K_n\circ\tau_x ((u\dvtx n)_{\alpha
} )$ on the cube $[(\mathbf{j}-\mathbf{1}),\mathbf{j}]$. This
point of view is the basic step in the approximation by the
Riemann--Stieltjes integral of $B(((\cdot)\dvtx \mathbf{n})_{\mathbf
{1}-\alpha})$ with respect to the function $(K_n\circ\tau_x)_{\alpha
}$ for each $\alpha\in\{0,1\}^d$.

\begin{Lemma} \label{lemma4}
\begin{eqnarray*}
\sup_{x\in[0,{h}^{-1}]^d}\bigl|Y_{n,2}^{(+)}(x)-Y_{n,3}^{(+)}(x)\bigr|=\mathrm{o}_P
\biggl(\frac{1}{\log(n)} \biggr) ,
\end{eqnarray*}
where the process $Y_{n,3}^{(+)}$ is defined by\vspace*{-1pt}
%
\begin{equation}
\label{yn3+} Y_{n,3}^{(+)}(x)\stackrel{\mathcal{D}} {=}
\frac{(2\pi)^d}{\|K\|
_{L^2}}\int_{[0,(a_nh)^{-1}]^d}K(x-u) \,\mathrm{d}B(u).\vspace*{-1pt}
\end{equation}
\end{Lemma}

We obtain from its proof in Section~\ref{sec7.2} that Lemma~\ref
{lemma4} holds under the condition $\frac{\log(n)}{nh^{d}}=\mathrm{o}(1)$,
which follows from the assumptions of Theorem~\ref{GWS}.
In the final step we show that the difference\vspace*{-1pt}
\[
Y^{(+)}(x)-Y^{(+)}_{n,3}(x)= \frac{(2\pi)^{d}}{\|K\|_{L^2}} \int
_{\R
^d_+}I_{\R^d_+\backslash[0,(a_nh)^{-1}]^d}(u) K(x-u) \,\mathrm{d}B(u)\vspace*{-1pt}
\]
is asymptotically negligible.

\begin{Lemma}\label{InfInt}\vspace*{-1pt}
\[
\sup_{x\in[0,h^{-1}]^d}\bigl|Y_{n,3}(x)-Y(x)\bigr|=
\mathrm{o}_P\bigl((\log(n)^{-1}\bigr).\vspace*{-1pt}
\]
\end{Lemma}

\subsection{Proof of Corollary \texorpdfstring{\protect\ref{cor1}}{1}}

The assertion follows from the estimate
%
\begin{eqnarray}
\label{e1} \sup_{[0,1]^d}\bigl\llvert f(x)-\E\hat
f_n(x)\bigr\rrvert =\mathrm{o} \bigl(h^{-\beta} \bigl({h^dn^da_n^d}
\bigr)^{-1/2} \bigr).\vspace*{-1pt}
\end{eqnarray}
To prove \eqref{e1} we use the representation \eqref{esti} and obtain
by a straightforward calculation\vspace*{-1pt}
\begin{eqnarray*}
\mathbb{E}\hat f_n(x)&=&\frac{1}{(2\pi)^{d}n^da_n^dh^{d}}\sum
_{\mathbf{k}\in G_{\mathbf{n}}}(f*\psi) (x_{\mathbf{k}}) \cdot K_n
\biggl((x-x_{\mathbf{k}})\frac{1}{h} \biggr)
\\[-1pt]
&=&\frac{1}{(2\pi)^{d}h^{d}}\int_{[-\sfrac{1}{a_n},\sfrac
{1}{a_n}]^d}(f*\psi) (z) \cdot
K_n \biggl((x-z)\frac{1}{h} \biggr) \,\mathrm{d}z +
R_{n,1}(x)
\\[-1pt]
&=&\frac{1}{(2\pi)^{d}}\int_{\R^d}(f*\psi) (z\cdot
h) \cdot K_n \biggl(\frac{x}{h}-z \biggr) \,\mathrm{d}z
+R_{n,1}(x)+R_{n,2}(x),
\end{eqnarray*}
where the term\vspace*{-1pt}
%
\begin{eqnarray}
\label{Rn1} R_{n,1}(x) &=& \frac{1}{(2\pi)^{d}h^d}\sum
_{\mathbf{k}\in
G_{\mathbf{n}-\mathbf{1}}}\int_{[x_{\mathbf{k}},x_{\mathbf{k+1}}]} \biggl\{ (f*\psi)
(x_{\mathbf{k}})K_n \biggl(\frac{x-x_{\mathbf
{k}}}{h} \biggr)\nonumber\\[-8pt]\\[-8pt]
&&\hphantom{\frac{1}{(2\pi)^{d}h^d}\sum
_{\mathbf{k}\in
G_{\mathbf{n}-\mathbf{1}}}\int_{[x_{\mathbf{k}},x_{\mathbf{k+1}}]} \biggl\{}{}-(f*\psi)
(z)K_{n} \biggl(\frac{x-z}{h} \biggr) \biggr\} \,\mathrm{d}z\nonumber
\end{eqnarray}
denotes the ``error'' in the integral approximation and
\begin{eqnarray*}
R_{n,2}(x):=\frac{1}{(2\pi)^{d}h^{d}}\int_{ ([-\sfrac
{1}{a_n},\sfrac{1}{a_n}]^d )^{C}}(f*\psi)
(z) K_n \biggl((x-z)\frac{1}{h} \biggr) \,\mathrm{d}z.
\end{eqnarray*}
An application of the Plancherel identity (see for example Folland \cite{folland1984}, Theorem~8.29) gives (observing
Assumption A1 and A3)
\begin{eqnarray*}
\mathbb{E}\hat f_n(x) & =&\frac{1}{(2\pi)^{\sfrac{d}{2}}h^{d}}\int
_{\R^d}\mathcal {F}f \bigl(h^{-1}\xi \bigr)
\mathcal{F}\psi \bigl(h^{-1}\xi \bigr) \frac{\mathcal{F} \eta(\xi) }{\mathcal{F} \psi(\sfrac{\xi}{h})} \exp
\bigl(\mathrm{i}h^{-1}x^T\xi \bigr) \,\mathrm{d}\xi\\
&&{}+R_{n,1}(x)+R_{n,2}(x)
\\
&=&\frac{1}{(2\pi)^{\sfrac{d}{2}}}\int_{\R^d}\mathcal{F}f(\xi)
\cdot\mathcal{F} \eta(\xi h)\exp\bigl(\mathrm{i}x^T\xi\bigr) \,\mathrm{d}\xi
+R_{n,1}(x)+R_{n,2}(x)
\\
&=& f(x) + R_{n,1}(x)+R_{n,2}(x) +
R_{n,3}(x)+R_{n,4}(x),
\end{eqnarray*}
where
\begin{eqnarray*}
R_{n,3}(x) &=& \frac{1}{(2\pi)^{\sfrac{d}{2}}}\int_{([-\sfrac
{D}{h},\sfrac{D}{h}]^d)^C}
\mathcal{F}f(\xi) \exp(\mathrm{i}x\xi) \,\mathrm{d}\xi,
\\
R_{n,4}(x) &=& \frac{1}{(2\pi)^{\sfrac{d}{2}}}\int_{[-\sfrac
{1}{h},\sfrac{1}{h}]^d\backslash[-\sfrac{D}{h},\sfrac
{D}{h}]^d}
\mathcal{F}f(\xi) \cdot\mathcal{F} \eta(\xi h)\exp(\mathrm{i}x\xi) \,\mathrm{d}\xi.
\end{eqnarray*}
We further obtain from Assumption \ref{assC}
\begin{eqnarray*}
\biggl|\int_{\{\xi_j>\sfrac{D}{h}\}}\mathcal{F}f(\xi) \exp(-\mathrm{i}x\xi) \,\mathrm{d}\xi \biggr|\leq
\frac{1}{D^{\gamma}}\int_{\{\xi_j>\sfrac{D}{h}\}}\bigl|\mathcal{F}f(\xi )\bigr|(h
\xi_j)^{\gamma} \,\mathrm{d}\xi =\mathrm{o}\bigl(h^{\gamma}\bigr), %
\end{eqnarray*}
and finally $|R_{n,3}(x)|
\leq\sum_{j=1}^d\int_{\{\xi_j>\sfrac{D}{h}\}}|\mathcal{F}f(\xi)|
\,\mathrm{d}\xi=\mathrm{o}(h^{\gamma})$.
With the same arguments it follows $R_{n,4}(x) =\mathrm{o}(h^{\gamma})$,
since $\llvert \mathcal{F}\eta(\xi h)\rrvert \leq1$ for all $\xi\in
\R^d$. Define $\mathcal{A}_n= ([-\frac{1}{a_n},\frac
{1}{a_n}]^d)^{C}$, then we
obtain from the representation \eqref{expan} the estimate
\begin{eqnarray*}
\bigl| R_{n,2}(x)\bigr| &\leq&\frac{1}{(2\pi)^{d}h^{\beta+d}} \biggl(\int_{\mathcal{A}_n}
\bigl|(f*\psi) (z)\bigr|^2 \,\mathrm{d}z \biggr)^{\sfrac{1}{2}} \\
&&{}\times\biggl[ \biggl(\int
_{ \mathcal{A}_n} \biggl|K \biggl((x-z)\frac{1}{h} \biggr)\biggr|^2
\,\mathrm{d}z \biggr)^{\sfrac{1}{2}}
\\
&&\hphantom{{}\times\biggl[}{} + \biggl(\int_{ \mathcal{A}_n }\biggl|\bigl(h^{\beta}K_n-K
\bigr) \biggl((x-z)\frac{1}{h} \biggr)\biggr|^2 \,\mathrm{d}z
\biggr)^{\sfrac{1}{2}} \biggr]
\\
&=&\frac{1}{(2\pi)^{d}h^{\beta}} \biggl(\int_{\mathcal{A}_n } \bigl|(f*
\psi) (z)\bigr|^2 \,\mathrm{d}z \biggr)^{\sfrac{1}{2}} \bigl(\mathrm{O} \bigl(h^da_n^{\sfrac
{d}{2}}
\bigr)+\mathrm{O} \bigl(h^{\mu_1+d}a_n^{\sfrac{d}{2}} \bigr) \bigr) = \mathrm{O}
\biggl(\frac{a_n^{\nu+\sfrac{d}{2}}}{h^{\beta}} \biggr)
\end{eqnarray*}
uniformly with respect to $x\in[0,1]^d$.
Note that by Assumption \ref{assC} we have $f\in\mathcal{W}^{ \lfloor m
\rfloor}(\R^d)$ and since
$m>2+\frac{d}{2}$ Sobolev's Embedding Theorem (Folland \cite{folland1984},
Theorem~8.54) implies\vadjust{\goodbreak} the existence of a function $\tilde f\in C^2(\R
^d)$ with $f=\tilde f$ almost everywhere. Observing that the
convolution function $\psi$ is integrable gives
$ \partial^{\alpha}(f*\psi)=(\partial^{\alpha}f)*\psi\in C(\R^d)$
for all $\alpha\in\{0,1,2\}$ with $|\alpha|\leq2$ (see for example
Folland \cite{folland1984}, Proposition~8.10), which
justifies the application
of Taylor's Theorem. Straightforward but tedious calculations give for
the remaining term \eqref{Rn1}
$
R_{n,1}(x)= \mathrm{O}(\frac{1}{n^3a_n^3h^{\beta+2}}) +\mathrm{O}(\frac{a_n^{\nu
}}{nh^{\beta}})
$
uniformly with respect to $x \in[0, 1]^d$.

\subsection{Proofs of Theorem \texorpdfstring{\protect\ref{GWSZeit}}{2} and Corollary
\texorpdfstring{\protect\ref{cor2}}{2}}

First we will show that the kernel $\check K_n$ satisfies conditions B1
and B2, with the kernel $\check K$
defined \eqref{kcheck}.
If Assumption \hyperref[assBB]{$\check \mathrm{B}$} holds we have
\begin{eqnarray*}
&&\int_{\R^d} \biggl(\frac{h^{\beta}}{\mathcal{F}_d\psi(\sfrac{\xi
}{h})}-\Psi(\xi) \biggr)
\mathcal{F}_{d+1}\check\eta(\xi,\tau)\exp \bigl(\mathrm{i}x^T\xi\bigr)
\,\mathrm{d}\xi
\\
&&\quad =\sum_{p=1}^{M-1}h^{\mu_p}\int
_{\R^d}\Psi_{p}(\xi)\mathcal {F}_{d+1}
\check\eta(\xi,\tau)\exp\bigl(\mathrm{i}x^T\xi\bigr) \,\mathrm{d}\xi\\
&&\qquad {} +h^{\mu_M}\int
_{\R^{d}}\Psi_{M,n}(\xi)\mathcal{F}_{d+1}
\check \eta(\xi,\tau)\exp\bigl(\mathrm{i}x^T\xi\bigr) \,\mathrm{d}\xi,
\end{eqnarray*}
which implies
\begin{eqnarray*}
h^{\beta} \check K_n(x,t)-\check K(x,t)&=&\sum
_{p=1}^{M-1}h^{\mu
_p}\int_{\R^{d+1}}
\Psi_{p}(\xi) \mathcal{F}_{d+1}\check\eta(\xi,\tau)\exp
\bigl(\mathrm{i}x^T\xi+\mathrm{i}t\tau\bigr) \,\mathrm{d}(\xi ,\tau)
\\
&&{}+h^{\mu_M}\int_{\R^{d+1}}
\Psi_{M,n}(\xi)\mathcal{F}_{d+1}\check \eta(\xi,\tau)\exp
\bigl(\mathrm{i}x^T\xi+\mathrm{i}t\tau\bigr)\, \mathrm{d}(\xi,\tau).
\end{eqnarray*}
A careful inspection of the proofs of Theorem~\ref{GWS} and
Corollary~\ref{cor1} shows that the arguments can be transferred
to the time-dependent case if the increase of $n$ and $m(n)$ as well as
the decrease of $a_n, b_{m(n)}$, $h$ and $h_t$ are
balanced as given in the assumptions of the theorem. The details are
omitted for the sake of brevity.

\section{Proof of auxiliary results} \label{sec7}
\setcounter{equation}{0}

\subsection{Proof Lemma \texorpdfstring{\protect\ref{Rio}}{2}} \label{sec7.1}

Define $S_{\mathbf{k}}:=\sum_{\mathbf{j}\in G^+_{\mathbf
{k}}}\varepsilon_{\mathbf{j}}$, set $S_{\mathbf{j}}\equiv0$
if $\min\{j_1,\ldots,j_d\}=0$ and recall the definition of
$Y_n^{(+)}$ and $\tau_x$ in \eqref{Yn+} and before Lemma~\ref{Rio},
respectively. In a first step we will replace the errors $\varepsilon
_{k}$ by increments given in terms of partial sums $S_{\mathbf
{k}-\alpha}$ for $\alpha\in\{0,1\}^d$. To be precise, we use the
representation
\begin{eqnarray*}
\varepsilon_{\mathbf{k}}=\sum_{\underline{\alpha}\in\{0,1\}
^d}(-1)^{|\underline{\alpha}|}S_{({\mathbf{k}}-\underline{\alpha})}
=\sum_{\underline{\alpha}\in\{0,1\}^d}(-1)^{|\underline{\alpha}|} S_{ ((\mathbf{k}-\mathbf{1}):\mathbf{k} )_{\underline
{\alpha}}}.
\end{eqnarray*}
A straightforward calculation gives
\begin{eqnarray*}
Y_n^{(+)}(x)&:=&\frac{h^{\beta}}{\sigma\|K\|_{2}\sqrt {n^dh^da_n^d}}\sum
_{\mathbf{k}\in G_n^+} K_n\circ\tau_x(\mathbf{k}-
\mathbf{1})\sum_{\alpha\in\{0,1\}
^d}(-1)^{|\alpha|}
S_{\mathbf{k}-\alpha}
\\
&=&\frac{h^{\beta}}{\sigma\|K\|_{2}\sqrt{n^dh^da_n^d}}\sum_{\alpha
\in\{0,1\}^d}(-1)^{|\alpha|}
\sum_{\mathbf{k}\in G_n^+} K_n\circ
\tau_x(\mathbf{k}-\mathbf{1}) S_{((\mathbf{k}-\mathbf{1}):\mathbf{k})\alpha}
\\
&=&\frac{h^{\beta}}{\sigma\|K\|_{2}\sqrt{n^dh^da_n^d}} \\
&&{}\times\biggl(\sum_{\alpha\in\{0,1\}^d}(-1)^{|\alpha|}
\sum_{\mathbf{k}\in G_n^+} \bigl(K_n\circ
\tau_x(\mathbf{k}-\mathbf{1})-K_n\circ
\tau_x\bigl(\bigl((\mathbf {k}-\mathbf{1})\dvtx \mathbf{k}
\bigr)_\alpha\bigr)\bigr) S_{((\mathbf{k}-\mathbf{1}):\mathbf{k})\alpha}
\\
&&\hphantom{{}\times\biggl(}{}+ \sum_{\alpha\in\{0,1\}^d}(-1)^{|\alpha|} \sum
_{\mathbf{k}\in G_n^+} K_n\circ\tau_x
\bigl(\bigl((\mathbf{k}-\mathbf{1})\dvtx \mathbf{k} \bigr)_\alpha\bigr)
S_{((\mathbf{k}-\mathbf{1}):\mathbf{k})\alpha} \biggr).
\end{eqnarray*}
Now we can make use of Proposition~6 and Proposition~3 of Owen \cite{owen2005} to rewrite the sums, such that the
increments given in terms
of partial sums can be expressed by increments given in terms of the
kernel $K_n$. We obtain
\begin{eqnarray*}
Y_n^{(+)}(x)&=&\frac{h^{\beta}}{\sigma\|K\|_{2}\sqrt {n^dh^da_n^d}} \\
&&{}\times\biggl[K_n
\circ\tau_x(\mathbf{n})S_{(\mathbf{n})}
\\
&&\hphantom{{}\times \biggl[}{} +\sum_{\alpha\in\{0,1\}^d}(-1)^{|\alpha|} \sum
_{\mathbf{k}\in G_n^+}\sum_{\beta\in\{0,1\}^d\backslash\{
\mathbf{0}\}}(-1)^{|\beta|}
\Delta_{\beta} \bigl(K_n\circ\tau_x;
\mathbf{k}-\mathbf{1}, \\
&&\hphantom{{}\times \biggl[{}+\sum_{\alpha\in\{0,1\}^d}(-1)^{|\alpha|} \sum
_{\mathbf{k}\in G_n^+}\sum_{\beta\in\{0,1\}^d\backslash\{
\mathbf{0}\}}(-1)^{|\beta|}
\Delta_{\beta} \bigl(}\bigl((\mathbf {k}-\mathbf{1})\dvtx \mathbf{k}
\bigr)_\alpha\bigr)S_{((\mathbf{k}-\mathbf{1}):\mathbf{k})\alpha} \biggr].
\end{eqnarray*}
The quantity
$\Delta_{\beta} (K_n\circ\tau_x;\mathbf{k}-\mathbf{1}, ((\mathbf
{k}-\mathbf{1})\dvtx \mathbf{k}
)_\alpha)$
can only take values different from zero if $\alpha\leq\mathbf
{1}-\beta$. Note that for $\alpha\leq\mathbf{1}-\beta$ the
equality $(\mathbf{k})_{\beta}=
((\mathbf{k}-\mathbf{1})\dvtx \mathbf{k})_{\alpha} )_{\beta}$ holds
which implies that in this case we also have $  [ (\mathbf
{k}-\mathbf{1} )_{\beta},  (((\mathbf{k}-\mathbf
{1})\dvtx \mathbf{k})_{\alpha}  )_{\beta}  ]=[(\mathbf
{k}-\mathbf{1})_{\beta},(\mathbf{k})_{\beta}]$. We further obtain
\begin{eqnarray*}
Y_n^{(+)}(x) &=&\frac{h^{\beta}}{\sigma\|K\|_{2}\sqrt{n^dh^da_n^d}} \\
&&{}\times\biggl[K_n
\circ\tau_x(\mathbf{n})S_{(\mathbf{n})} \\
&&\hphantom{{}\times\biggl[}{}+\sum
_{\beta\in\{0,1\}^d\backslash\{\mathbf{0}\}}(-1)^{|\beta|} \sum_{\mathbf{k}\in G_{\mathbf{n}}^+}
\sum_{\tilde{\alpha}\in\{
0,1\}^{d-|\beta|}}(-1)^{|\tilde{\alpha}|}
\\
&&\hphantom{{}\times\biggl[{}+\sum
_{\beta\in\{0,1\}^d\backslash\{\mathbf{0}\}}(-1)^{|\beta|} \sum_{\mathbf{k}\in G_{\mathbf{n}}^+}
\sum_{\tilde{\alpha}\in\{
0,1\}^{d-|\beta|}}}{} \times\Delta_{\beta} \bigl( K_n\circ
\tau_x ; \mathbf{k}-\mathbf{1} , (\mathbf{k})_{\beta}\dvtx\\
&&\hphantom{{}\times\biggl[{}+\sum
_{\beta\in\{0,1\}^d\backslash\{\mathbf{0}\}}(-1)^{|\beta|} \sum_{\mathbf{k}\in G_{\mathbf{n}}^+}
\sum_{\tilde{\alpha}\in\{
0,1\}^{d-|\beta|}}{} \times\Delta_{\beta} \bigl(}{} \bigl((
\mathbf{k}-\mathbf {1})_{\mathbf{1}-\beta} \dvtx (\mathbf{k})_{\mathbf{1}-\beta}
\bigr)_{\tilde{\alpha}} \bigr) \\
&&\hphantom{{}\times\biggl[{}+\sum
_{\beta\in\{0,1\}^d\backslash\{\mathbf{0}\}}(-1)^{|\beta|} \sum_{\mathbf{k}\in G_{\mathbf{n}}^+}
\sum_{\tilde{\alpha}\in\{
0,1\}^{d-|\beta|}}}{} \times S_{(\mathbf{k})_{\beta}:((\mathbf{k}-\mathbf{1})_{\mathbf{1}-\beta}
:(\mathbf{k})_{\mathbf{1}-\beta})_{\tilde{\alpha}} } \biggr] .
\end{eqnarray*}
The alternating sum with respect to the index
$\tilde{\alpha}$ can be written as an increment $\Delta$ as defined
in \eqref{Delta} which then defines a telescope sum according to
Owen \cite{owen2005}, Proposition~2. Taking into
account that
$S(\mathbf{k})\equiv0$ if $k_j=0$ for at least one $j\in\{1,\ldots
,d\}$ gives
\begin{eqnarray*}
Y_n^{(+)}(x)=\frac{h^{\beta}}{\sigma\|K\|_{2}\sqrt{n^dh^da_n^d}} \sum
_{\beta\in\{0,1\}^d}(-1)^{|\beta|} \sum_{\mathbf{j}\in G_{\mathbf{n}}^{|\beta|,+}}
\Delta_{\beta} ( K_n\circ\tau_x
;I_{\mathbf{j}} %
)\cdot S_{\mathbf{j}:(\mathbf{n})_{\mathbf{1}-\beta}}.
\end{eqnarray*}
With the definitions
$X(A):=\sum_{\mathbf{k}\in A\subset\mathbb{Z}^d}X_{\mathbf{k}}$
for any subset $A\in\mathbb{Z}^d$, we can rewrite these partial sums
as set-indexed partial sums with index class
$n\cdot\mathscr{S}$, where $\mathscr{S}:=\{(0,\gamma] | 0<\gamma
_j\leq1, 1\leq j\leq d\}$ and $n\cdot\mathscr{S}:=\{n\cdot S | S\in
\mathscr{S}\}$. It follows
directly that $\mathscr{S}$ is a sufficiently smooth VC-class of sets,
which justifies the application of Theorem~1 in Rio
\cite{rio1993}.
Therefore there exists a version of a Brownian sheet on $[0,\infty
)^d$, say
$B_1$, such that
%
\begin{eqnarray}
\label{rio} \sup_{\mathbf{k}\in G_{\mathbf{n}}^+}\biggl\llvert \frac{S_{\mathbf
{k}}}{\sigma}-B_1(
\mathbf{k})\biggr\rrvert =\mathrm{O} \bigl(\bigl(\log(n)\bigr)^{\sfrac
{1}{2}}n^{\vfrac{d-\delta}{2}}
\bigr)\qquad  \mbox{a.s.}
\end{eqnarray}
Recalling the definition of $I_{\mathbf{j}}$ in \eqref{Interval} we
further obtain
\begin{eqnarray*}
&&Y_{n}^{(+)}(x)-Y_{n,1}^{(+)}(x)\\
&&\quad =
\frac{h^{\beta}}{\|K\|_{2}\sqrt {n^dh^da_n^d}} \\
&&\qquad {}\times\sum_{\beta\in\{0,1\}^d}(-1)^{|\beta|} \sum_{\mathbf{j}\in G_{\mathbf{n}}^{|\beta|,+}} \Delta_{\beta} (
K_n\circ\tau_x ;I_{\mathbf{j}} %
)\cdot
\biggl( \frac{1}{\sigma}S_{\mathbf{j}: (\mathbf
{n})_{\mathbf{1}-\beta}}-B_1 \bigl(\mathbf{j}\dvtx (
\mathbf{n})_{\mathbf
{1}-\beta} \bigr) \biggr).
\end{eqnarray*}
The estimate \eqref{rio} implies the existence of a constant $C\in\R
_+$ such that
\begin{eqnarray*}
&&\bigl|Y_{n}^{(+)}(x)-Y_{n,1}^{(+)}(x)\bigr|\\
&&\quad \leq C
\cdot\sqrt{\frac{\log
(n)}{n^{\delta}h^{\delta}a_n^{\delta}}}h^{\beta} \biggl[\sum
_{\gamma\in\{0,1\}^d, |\gamma|=1}\int_{[0,(a_nh)^{-1}]^d}(u)_{\gamma}^{\vfrac{d-\delta}{2}}\bigl|
\partial ^{\mathbf{1}}K_n(x-u)\bigr|\, \mathrm{d}u
\\
&&\hphantom{\quad \leq C
\cdot\sqrt{\frac{\log
(n)}{n^{\delta}h^{\delta}a_n^{\delta}}}h^{\beta} \biggl[}{} +\sum_{\beta\in\{0,1\}^d\backslash\{\mathbf{0},\mathbf{1}\}
}\int_{[0,(a_nh)^{-1}]^{|\beta|}}
\bigl\llvert \partial^{\beta}K_n \bigl( \bigl(x-
\bigl(u\dvtx (a_nh)^{-1}\mathbf{1}\bigr) \bigr)_{\beta}
\bigr) \bigr\rrvert (\mathrm{d}u)_{\beta}
\\
&&\hphantom{\quad \leq C
\cdot\sqrt{\frac{\log
(n)}{n^{\delta}h^{\delta}a_n^{\delta}}}h^{\beta} \biggl[}{} +\bigl|K_n\bigl(x-(a_nh)^{-1}
\mathbf{1}\bigr)\bigr| \biggr] \qquad \mbox{a.s.}
\end{eqnarray*}
It follows from Assumption \ref{assB} that the function $u\mapsto(u)_{\gamma
}^{\sfrac{|\alpha|}{2}}\,\partial^{\alpha}K(u)$ is integrable on
$\mathbb{R}^d$ for all $\alpha\in\{0,1\}^d$ such that\vspace*{-1.7pt}
\begin{eqnarray*}
\int_{[0,(a_nh)^{-1}]^d}(u)_{\gamma}^{\vfrac{d-\delta}{2}}\bigl|\partial
^{\mathbf{1}}K_n(x-u)\bigr| \,\mathrm{d}u=\mathrm{O}\bigl(h^{\vfrac{\delta-d}{2}-\beta}\bigr)
\end{eqnarray*}
and\vspace*{-1.7pt}
\begin{eqnarray*}
\int_{[0,(a_nh)^{-1}]^{|\beta|}} \bigl| \partial^{\beta}K_n \bigl(
\bigl(x-\bigl(u\dvtx (a_nh)^{-1}\mathbf{1}\bigr)
\bigr)_{\beta} \bigr) \bigr| (\mathrm{d}u)_{\beta}+\bigl|K_n
\bigl(x-(a_nh)^{-1}\mathbf{1}\bigr)\bigr| =\mathrm{O}
\bigl((a_nh)^{\sfrac{d}{2}}h^{-\beta}\bigr).
\end{eqnarray*}

Note that for sufficiently large $n$ such that $a_n<\frac{1}{2}$ we obtain
$
-\frac{1}{2a_nh}\geq x_j-(a_nh)^{-1}=\frac{a_n-1}{a_nh}
$
uniformly with respect to $j$ (note that $x_j\in[0,h^{-1}]$).
Let $\tilde B$ be a continuous version of $B_1$. We set $\tilde
B(t)\equiv0$ if $t_j<0$ for at least one index $j \in\{1,\ldots,d\}$
and let $\{\tilde B_{\alpha} | \alpha\in\{0,1\}^d\}$ be $2^d$
mutually independent copies of $\tilde B$. For $t\in\R^d$ define\vspace*{-1.7pt}
\begin{eqnarray*}
B_{\alpha}(t):=\tilde B_{\alpha}\bigl((-1)^{\alpha_1}t_1,(-1)^{\alpha
_2}t_2,
\ldots,(-1)^{\alpha_d}t_d\bigr),
\end{eqnarray*}
then the process
$ \{ B(t):=\sum_{\alpha\in\{0,1\}^d}B_{\alpha}(t) | t\in\R^d\}
$
is a Wiener field on $\R^d$.

\subsection{Proof of Lemma \texorpdfstring{\protect\ref{lemma4}}{4}}\label{sec7.2}
Note that $\partial^{\alpha}K$ exists and is integrable for each
$\alpha\in\{0,1\}^d$. Consequently, the
kernel $K$ is of bounded variation on $[{0},(a_nh)^{-1}]^d$ in the
sense of Hardy Krause for each fixed $n$
(see Owen \cite{owen2005}, Definition~2). Therefore
an application of
integration by parts for the Wiener integral
(note that the kernel $K$ has not necessarily a compact support) and
rescaling of the Brownian sheet $Y_{n,3}^{(+)}$ yields\vspace*{-1.7pt}
\begin{eqnarray*}
Y_{n,3}^{(+)}(x)&\stackrel{\mathcal{D}} {=}&\sum
_{\alpha\in\{0,1\}
^d\backslash\{\mathbf{0}\}}(-1)^{|\alpha|} \int_{[0,(a_nh)^{-1}]^{|\alpha|}}B
\bigl( \bigl(u\dvtx (a_nh)^{-1}\mathbf {1}
\bigr)_{\alpha} \bigr) \,\mathrm{d}K \bigl(x- \bigl(u\dvtx (a_nh)^{-1}
\mathbf{1} \bigr)_{\alpha} \bigr)
\\[-1.7pt]
&&{}+\Delta \bigl(K(x-\cdot)\cdot B(\cdot),\bigl[0,(a_nh)^{-1}
\bigr]^d \bigr)
\\[-1.7pt]
&=&\sum_{\alpha\in\{0,1\}^d\backslash\{\mathbf{0}\}}(-1)^{|\alpha|} \int
_{[0,(a_nh)^{-1}]^{|\alpha|}}B \bigl( \bigl(u\dvtx (a_nh)^{-1}
\mathbf {1} \bigr)_{\alpha} \bigr) \,\partial^{\alpha}K \bigl(x-
\bigl(u\dvtx (a_nh)^{-1}\mathbf{1} \bigr)_{\alpha} \bigr)
(\mathrm{d}u)_{\alpha}
\\[-1.7pt]
&&{}+\Delta \bigl(K(x-\cdot)\cdot B(\cdot),\bigl[0,(a_nh)^{-1}
\bigr]^d \bigr).
\end{eqnarray*}
Recalling the definition of $Y_{n,2}^{(+)}(x)$ and identity (15) in
Owen \cite{owen2005} we can replace the
increments by the corresponding integrals, that is\vspace*{-1.7pt}
\begin{eqnarray*}
Y_{n,2}^{(+)}(x) &\stackrel{\mathcal{D}} {=}&\sum
_{\alpha\in\{0,1\}^d\backslash\{
\mathbf{0}\}}(-1)^{|\alpha|}\\[-1.7pt]
&&{}\times\sum_{\mathbf{k}\in
G_{\mathbf{n}-\mathbf{1}}^{|\alpha|,+}}
\int_{[(na_nh)^{-1}(\mathbf{k}-\mathbf{1})_{\alpha
},(na_nh)^{-1}(\mathbf{k})_{\alpha}]} \partial^{\alpha}K \bigl(x-
\bigl(u\dvtx (a_nh)^{-1}\mathbf{1} \bigr)_{\alpha} \bigr)
(\mathrm{d}u)_{\alpha}
\\
&&\hphantom{{}\times\sum_{\mathbf{k}\in
G_{\mathbf{n}-\mathbf{1}}^{|\alpha|,+}}
\int_{[(na_nh)^{-1}(\mathbf{k}-\mathbf{1})_{\alpha
},(na_nh)^{-1}(\mathbf{k})_{\alpha}]}}{} \times B \bigl( \bigl((na_nh)^{-1}
\mathbf{k}\dvtx (a_nh)^{-1}\mathbf {1} \bigr)_{\alpha}
\bigr) \\
&&\hphantom{{}\times\sum_{\mathbf{k}\in
G_{\mathbf{n}-\mathbf{1}}^{|\alpha|,+}}
\int_{[(na_nh)^{-1}(\mathbf{k}-\mathbf{1})_{\alpha
},(na_nh)^{-1}(\mathbf{k})_{\alpha}]}}{}+\Delta \bigl(K(x-\cdot)\cdot B(\cdot),\bigl[0,(a_nh)^{-1}
\bigr]^d \bigr)
\\
&=&Y^+_{n,3}(x)+R_{n,\mathrm{SI}}(x),
\end{eqnarray*}
where the remainder $R_{n,\mathrm{SI}}(x)$ is defined in an obvious manner.
From the modulus of continuity for the Brownian Sheet (see Khoshnevisan \cite{khosh2002}, Theorem~3.2.1) it follows that for
$a,b\in\R^d$\vspace*{1.5pt}
%
\begin{eqnarray}
\label{modcont} \limsup_{\delta\rightarrow0+}\sup_{s,t\in[a,b], \|s-t\|_{\infty
}<\delta}
\frac{|B(s)-B(t)|}{\sqrt{\delta\log(\sfrac{1}{\delta})}} \leq24 \cdot d \|b\|_{\infty}^{d/2},
\end{eqnarray}
which yields\vspace*{1.5pt}
\begin{eqnarray*}
\bigl|Y_{n,2}^{(+)}(x)-Y_{n,3}^{(+)}(x)\bigr|&=&\bigl|R_{n,\mathrm{SI}}(x)\bigr|\\
&\leq&\sup_{\delta
<\sfrac{1}{n}} \sup_{s,t\in [0,2 ]^d: \|s-t\|_{\infty
}\leq\delta}\bigl|B(s)-B(t)\bigr|
\\
&&{}\times\sqrt{\frac{\log(n)}{n}} \biggl[\int_{[0,(a-nh)^{-1}]^d}
\bigl((u)^{\mathbf{1}} \bigr)^{\sfrac
{1}{2}}\bigl|\partial^{\mathbf{1}}K(x-u)\bigr| \,\mathrm{d}u
+ \mathrm{O}\bigl(h^{-\vfrac{d-1}{2}}\bigr) \biggr]
\end{eqnarray*}
(note that the dominating term in $R_{n,\mathrm{SI}}(x)$ is given by the summand
where $|\alpha|=d$).
With the same arguments as in the proof of Lemma~\ref{Rio} we finally obtain\vspace*{1.5pt}
\begin{eqnarray*}
\bigl|{Y}_{n,2}^{(+)}(x)-Y_{n,3}^{(+)}(x)\bigr|
=\mathrm{O}_P \biggl(\sqrt{\frac{\ln(na_nh)}{nh^{d}}} \biggr),
\end{eqnarray*}
where we used the estimate \eqref{modcont} for the modulus of
continuity of the Brownian sheet (note that this estimate is
independent of $x$).

\subsection{Proof of Lemma \texorpdfstring{\protect\ref{InfInt}}{5}}
Integration by parts gives
%
\begin{eqnarray}
\label{delta} \Delta_{n,3}&:=&\bigl|Y_{n,3}^{(+)}(x)-Y^{(+)}(x)\bigr|
\nonumber\\
&\leq& \biggl|\int_{[0,\infty)^d\backslash[0,\afrac
{1}{a_nh}]^d}B(u)
\partial^{\mathbf{1}}K(x-u) \,\mathrm{d}u \biggr|
\nonumber\\
&&{}+ \biggl|\sum_{\alpha\in\{0,1\}^d\backslash\{\mathbf{0},\mathbf
{1}\}}(-1)^{|\alpha|}\nonumber\\[-10pt]\\[-6pt]
&&\hphantom{{}+ \biggl|\sum_{\alpha\in\{0,1\}^d\backslash\{\mathbf{0},\mathbf
{1}\}}}{}\times\int_{[0,\afrac{1}{a_nh}]^{|\alpha|}} B \bigl( \bigl(u\dvtx (a_nh)^{-1}
\mathbf{1} \bigr)_{\alpha} \bigr) \,\partial^{\alpha}K \bigl(x-
\bigl(u\dvtx (a_nh)^{-1}\mathbf{1} \bigr)_{\alpha} \bigr)
(\mathrm{d}u)_{\alpha} \biggr|
\nonumber\\
&&{}+\bigl\llvert \Delta \bigl(K(x-\cdot)B(\cdot);
\bigl[0,(a_nh)^{-1}\bigr]^d \bigr)\bigr\rrvert\nonumber\\
&:=&\bigl|\Delta_{n,3}^{(1)}(x)\bigr|+\bigl|\Delta_{n,3}^{(2)}(x)\bigr|+\bigl|
\Delta_{n,3}^{(3)}(x)\bigr|,\nonumber
\end{eqnarray}
where the processes $\Delta_{n,3}^{(j)}(x), j=1,2,3$ are defined in an
obvious manner. Let $n$ be sufficiently large such that $\frac
{1}{a_nh}\geq1$ and $a_n<\frac{1}{2}$.
Since $B(u)=0$ if $t_j=0$ for at least one index $j\in\{1,\ldots,d\}$
we have
\begin{eqnarray*}
\bigl|\Delta_{n,3}^{(3)}(x)\bigr|&=&\bigl\llvert K
\bigl(x-(a_nh)^{-1}\mathbf{1} \bigr)\cdot B
\bigl((a_nh)^{-1}\mathbf{1} \bigr)\bigr\rrvert
\\
&=&\sqrt{2d (a_nh )^{-d}\ln\bigl(d\ln
\bigl((a_nh)^{-1}\bigr)\bigr)} \frac{\llvert K (x-(a_nh)^{-1}\mathbf{1} )\rrvert |B
((a_nh)^{-1}\mathbf{1} )|}{\sqrt{2d (a_nh )^{-d}\ln
(d\ln((a_nh)^{-1}))}}.
\end{eqnarray*}
An application of the version of a law of the iterated logarithm given
in Theorem~3 of Paranjape and Park \cite{parpar1973} yields the estimate
\begin{eqnarray*}
\sup_{x\in[0,h^{-1}]}\bigl|\Delta_{n,3}^{(3)}(x)\bigr|&=&\mathrm{O}(1)
\cdot\sqrt {2d (a_nh )^{-d}\ln\bigl(d\ln
\bigl((a_nh)^{-1}\bigr)\bigr)}\sup_{x\in
[0,h^{-1}]}
\bigl\llvert K \bigl(x-(a_nh)^{-1}\mathbf{1} \bigr)\bigr
\rrvert
\\
&\leq& \mathrm{O}(1)\cdot\sqrt{2d (a_nh )^{-d}\ln\bigl(d\ln
\bigl((a_nh)^{-1}\bigr)\bigr)}\sup_{v\leq\afrac{a_n-1}{a_nh}}
\bigl\llvert K (v )\bigr\rrvert\\
& =&\mathrm{o} \biggl(\frac{1}{\log(n)} \biggr)\qquad  \mbox{a.s.}
\end{eqnarray*}
uniformly with respect to $x$.

To show that $\Delta_{n,3}^{(2)}(x)$ and $\Delta_{n,3}^{(1)}(x)$ are
asymptotically negligible we also apply the LIL for the Brownian
sheet.
For each summand, say $\Delta^{(2)}_{n,3,\alpha}$, in $\Delta
_{n,3}^{(2)}(x)$ $(|\alpha|<d)$ we have
\begin{eqnarray*}
\Delta_{n,3,\alpha}^{(2)}(x)&:=& \biggl| \int_{[0,\afrac
{1}{a_nh}]^{|\alpha|}} B
\bigl( \bigl(u\dvtx (a_nh)^{-1}\mathbf{1} \bigr)_{\alpha}
\bigr)\, \partial^{\alpha}K \bigl(x- \bigl(u\dvtx (a_nh)^{-1}
\mathbf{1} \bigr)_{\alpha} \bigr) (\mathrm{d}u)_{\alpha} \biggr|
\\
&=&(a_nh)^{-|\alpha|} \biggl| \int_{[0,1]^{|\alpha|}}
B \bigl( (u\dvtx \mathbf{1} )_{\alpha}(a_nh)^{-1} \bigr)
\,\partial^{\alpha}K \bigl(x- (u\dvtx \mathbf{1} )_{\alpha}(a_nh)^{-1}
\bigr) (\mathrm{d}u)_{\alpha} \biggr|. %
\end{eqnarray*}
Scaling of the Brownian sheet yields
\begin{eqnarray*}
\Delta_{n,3,\alpha}^{(2)}(x) &\stackrel{\mathcal{D}}
{=}&(a_nh)^{-\vfrac{2|\alpha|+d}{2}} \biggl| \int_{[0,1]^{|\alpha|}} B \bigl(
(u\dvtx \mathbf{1} )_{\alpha} \bigr) \,\partial^{\alpha}K \bigl(x- (u\dvtx
\mathbf{1} )_{\alpha}(a_nh)^{-1} \bigr)
(\mathrm{d}u)_{\alpha} \biggr|
\\
&=&\mathrm{O} \bigl( (a_nh)^{-\sfrac{d}{2}} \bigr) \biggl| \int
_{[0,\afrac
{1}{a_nh}]^{|\alpha|}} \partial^{\alpha}K \bigl(x- (u\dvtx \mathbf{1}
)_{\alpha}(a_nh)^{-1} \bigr) (\mathrm{d}u)_{\alpha} \biggr|
\qquad \mbox{a.s.}
\end{eqnarray*}
With the same arguments as in the proof of Lemma~\ref{Rio} we conclude
that the leading contributions are given by the quantities
$\Delta_{n,3,\alpha}^{(2)}(x)$, where $|\alpha|=d-1$. For $\alpha
=(0,1,\ldots,1) $ obtain
\begin{eqnarray*}
\sup_{x\in[0,h^{-1}]^d} \bigl|\Delta_{n,3,\alpha}^{(2)}(x) \bigr|=
\mathrm{O}_P \bigl( (a_nh)^{-\sfrac{d}{2}} \bigr) \sup
_{v\leq-\afrac{1}{2a_nh}}\int_{\R^{d-1}}\bigl\llvert \partial
^{\alpha}K (v,u_2,\ldots,u_d )\bigr\rrvert
\,\mathrm{d}(u_2,\ldots,u_d).
\end{eqnarray*}
This gives
$\sup_{x\in[0,h^{-1}]^d} |\Delta_{n,3,(0,1,\ldots
,1)}^{(2)}(x) |=\mathrm{o}(\frac{1}{\log(n)})$.
Applying the same argument to the other terms yields $\Delta
_{n,3}^{(2)}(x)=\mathrm{o}_P(\frac{1}{\log(n)})$ uniformly with respect to $x
\in[ 0,1/h]^d$. Finally, a similar argument gives for the remaining
term in \eqref{delta} $\Delta_{n,3}^{(1)}(x)=\mathrm{o}_P(\frac{1}{\log
(n)})$, which completes the proof of Lemma~\ref{InfInt}.

\subsection{Proof of Lemma \texorpdfstring{\protect\ref{K_nK}}{3}}
Note that we have
\begin{eqnarray*}
Y_{n,1}^{(+)}(x)\stackrel{\mathcal{D}} {=}\frac{h^{\beta}}{\sqrt {n^da_n^dh^d}\|K\|_{L^2(\R^d)}}
\sum_{\alpha\in\{0,1\}^d}(-1)^{|\alpha|}\sum
_{\mathbf{j}\in
G_n^{|\alpha|,+}}\Delta_{\alpha} ( \overline{\mathcal{F}f_{p}}
\circ\tau_x ;I_{\mathbf{j}} )B\bigl(\mathbf {j}\dvtx (
\mathbf{n})_{\mathbf{1}-\alpha}\bigr).
\end{eqnarray*}
The representation \eqref{expan} and the definition \eqref{Interval} yield
\begin{eqnarray*}
&&\bigl|Y_{n,1}^{(+)}(x)-Y_{n,2}^{+}(x)\bigr|
\\
&&\quad =\sum_{p=1}^{M-1}
\frac{h^{\mu_p}}{\sqrt{n^da_n^dh^d}} \biggl|\sum_{\alpha\in\{0,1\}^d}(-1)^{|\alpha|}
\sum_{\mathbf{j}\in G_{\mathbf{n}}^{|\alpha|,+}} \Delta_{\alpha} ( \overline{
\mathcal{F}f_{p}}\circ\tau_x ; I_{\mathbf{j}} )B
\bigl(\mathbf {j}\dvtx (\mathbf{n})_{\mathbf{1}-\alpha}\bigr) \biggr| +\mathrm{o}_P \biggl(
\frac{1}{\log(n)} \biggr).
\end{eqnarray*}
For each fixed $p$ we can now perform the approximation steps of the
previous lemmas and obtain
\begin{eqnarray*}
&& \log(n) \sup_{x\in[0,h^{-1}]^d} \biggl|\frac{1}{\sqrt{n^da_n^dh^d}}\sum
_{\alpha\in\{0,1\}
^d}(-1)^{|\alpha|} \sum_{\mathbf{j}\in G_{\mathbf{n}}^{|\alpha|,+}}
\Delta_{\alpha} ( \overline{\mathcal{F}f_{p}}\circ
\tau_x ; I_{\mathbf{j}} )B\bigl(\mathbf {j}\dvtx (
\mathbf{n})_{\mathbf{1}-\alpha}\bigr)
\\
&& \hphantom{\log(n) \sup_{x\in[0,h^{-1}]^d} |}{}-\int_{\R_+^d}\overline{\mathcal{F}f_p}(x-u)
\,\mathrm{d}B(u) \biggr|=\mathrm{o}_P (1 ).
\end{eqnarray*}
It can easily be shown that for all $p=1,\ldots,M-1$
\begin{eqnarray*}
\lim_{n\rightarrow\infty}n^da_n^dh^d
\Var \biggl(\frac
{1}{n^da_n^dh^d}\sum_{\mathbf{k}\in
G_{\mathbf{n}}}Y_{\mathbf{k}}
\overline{\mathcal{F}f_{p}} \biggl((x-x_{\mathbf{k}})
\frac{1}{h} \biggr) \biggr)=\sigma^2\|f_p
\|_{2}^2,
\end{eqnarray*}
where the limit does not depend on $x$.
We finally obtain, repeating the approximation steps given in the
previous lemmas for each of the $2^{d}-1$ remaining orthants
\begin{eqnarray*}
&&\hspace*{-5pt}\log(n) \sup_{x\in[0,h^{-1}]^d} \biggl|\frac{1}{\sqrt{n^da_n^dh^d}}\sum
_{\alpha,\gamma\in\{0,1\}
^d}(-1)^{|\alpha|} \sum_{\mathbf{j}\in G_{\mathbf{n}}^{|\alpha|,+}}
\Delta_{\alpha} \bigl( \overline{\mathcal{F}f_{p}}\circ
\tau_x ; (-1)^{\gamma}I_{\mathbf
{j}} \bigr)B
\bigl((-1)^{\gamma}\mathbf{j}\dvtx (\mathbf{n})_{\mathbf{1}-\alpha}\bigr)
\\
&&\hspace*{-5pt}\hphantom{\log(n) \sup_{x\in[0,h^{-1}]^d} |}{} -\int_{\R^d}\overline{\mathcal{F}f_p}(x-u) \,\mathrm{d}
B(u) \biggr|=\mathrm{o}_P(1).
\end{eqnarray*}
Note that
\begin{eqnarray*}
r(x-z) :=\mathbb{E} \biggl(\int_{\R^d}\overline{\mathcal
{F}f_p}(x-u) \,\mathrm{d} B(u)\int_{\R^d}
\mathcal{F}f_p(z-u) \,\mathrm{d} B(u) \biggr) =\int_{\R^d}f_p(x-z+u)f_p(u)
\,\mathrm{d}u
\end{eqnarray*}
and $\|r\|_1\leq\|f_p\|^2<\infty$.
The system of sets $\{[0,h^{-1}]^d | n\in\N\}$ is a blowing up system
of sets in the sense of Definition 14.1 in Piterbarg
\cite{piterbarg1996}.
If we define
\begin{eqnarray*}
Z_p(x)=\frac{1}{\|f_p\|_{L^2(\R^d)}}\int_{\R^d}\overline{
\mathcal {F}f_p}(x-u) \,\mathrm{d} B(u),
\end{eqnarray*}
then Theorem~1 in Bickel and Rosenblatt \cite{bicros1973b} gives
the asymptotic independence
of the scaled minimum and maximum of the process $Z_p$, which, with the
observation that $Z_p$ and $-Z_p$ have the same distribution
and an application of Theorem~14.1 in Piterbarg \cite
{piterbarg1996} yields that
for $G\sim\operatorname{Gumbel}(\ln(2),1)$
\begin{eqnarray*}
\sup_{x\in[0,h^{-1}]^d} \bigl( \bigl(\bigl|Z_p(x)\bigr|-\tilde
C_{n,3} \bigr)\tilde C_{n,3} \bigr)\stackrel{\mathcal{D}} {
\rightarrow}G \qquad \mbox{for } n\rightarrow\infty,
\end{eqnarray*}
where the constants $\tilde C_1, \tilde C_{n,2}$ and $\tilde C_{n,3}$
are given by
\begin{eqnarray*}
\tilde C_1&=&\operatorname{\mathbf{det}} \biggl( \biggl[\frac{1}{\|f_p\|_{L^2(\R
^d)}^2}\int
_{\R^d}\bigl|f_p(v)\bigr|^2v_iv_j
\,\mathrm{d}v \biggr] , i,j=1,\ldots,d \biggr),
\\
\tilde C_{n,2}&=&\sqrt{\frac{\tilde C_1}{(2\pi)^{d+1}}}
\frac{1}{h^d},
\\
\tilde C_{n,3}&=&\sqrt{2\ln(\tilde C_{n,2})}+
\frac{(d-1)\ln
(2\ln(\tilde C_{n,2}) )}{2\sqrt{2\ln(\tilde C_{n,2})}}.
\end{eqnarray*}

Since $h^{\mu_p}=\mathrm{o} (\frac{1}{\log n} )$ we obtain
$ h^{\mu_p}\sup_{x\in[0,h^{-1}]}|Z_p(x)|=\mathrm{o}_P ((\log n)^{-1/2})$ for
each $p=1,\ldots,M-1$,
which justifies the replacement of $h^{\beta}K_n$ by $K$.
Since the outer sum does not depend on $n$ this gives the desired result.

\subsection{Proof of Lemma \texorpdfstring{\protect\ref{Rand}}{1}}
With the same arguments as in the proof of the previous lemmas we can
replace the errors by combinations of partial sums and perform the same
approximation steps. In each replacement we
obtain at most a $d-1$-fold sum which yields the desired result right away.

\section*{Acknowledgements}
The authors would like to thank the Editor, the Associate Editor and
two referees for helpful
comments. Moreover, we thank Martina Stein, who typed parts of this
manuscript with considerable technical expertise.
This work has been supported in part by the
Collaborative Research Center ``Statistical modeling of nonlinear
dynamic processes'' (SFB 823,
Teilprojekt A1 and C4) of the German Research Foundation (DFG).



\printhistory


\begin{thebibliography}{32}


\bibitem{adorf1995}
\begin{barticle}[auto:STB|2014/01/06|10:16:28]
\bauthor{\bsnm{Adorf},~\bfnm{H.}\binits{H.}}
(\byear{1995}).
\btitle{Hubble space telescope image restauration in its fourth year}.
\bjournal{Inverse Problems}
\bvolume{11}
\bpages{639--653}.
\end{barticle}
\bptok{imsref}%
\endbibitem

\bibitem{berbocdes2009}
\begin{barticle}[mr]
\bauthor{\bsnm{Bertero},~\bfnm{M.}\binits{M.}},
\bauthor{\bsnm{Boccacci},~\bfnm{P.}\binits{P.}},
\bauthor{\bsnm{Desider{\`a}},~\bfnm{G.}\binits{G.}} \AND
\bauthor{\bsnm{Vicidomini},~\bfnm{G.}\binits{G.}}
(\byear{2009}).
\btitle{Image deblurring with {P}oisson data: From cells to galaxies}.
\bjournal{Inverse Problems}
\bvolume{25}
\bpages{123006, 26}.
\bid{doi={10.1088/0266-5611/25/12/123006}, issn={0266-5611}, mr={2565572}}
\end{barticle}
\bptok{imsref}%
\endbibitem

\bibitem{bicros1973b}
\begin{bincollection}[mr]
\bauthor{\bsnm{Bickel},~\bfnm{P.}\binits{P.}} \AND
\bauthor{\bsnm{Rosenblatt},~\bfnm{M.}\binits{M.}}
(\byear{1973}).
\btitle{Two-dimensional random fields}.
In \bbooktitle{Multivariate Analysis, {III} ({P}roc. {T}hird {I}nternat. {S}ympos., {W}right {S}tate {U}niv., {D}ayton, {O}hio, 1972)}
\bpages{3--15}.
\blocation{New York}:
\bpublisher{Academic Press}.
\bid{mr={0348832}}
\end{bincollection}
\bptok{imsref}%
\endbibitem

\bibitem{bicros1973a}
\begin{barticle}[mr]
\bauthor{\bsnm{Bickel},~\bfnm{P.~J.}\binits{P.J.}} \AND
\bauthor{\bsnm{Rosenblatt},~\bfnm{M.}\binits{M.}}
(\byear{1973}).
\btitle{On some global measures of the deviations of density function estimates}.
\bjournal{Ann. Statist.}
\bvolume{1}
\bpages{1071--1095}.
\bid{issn={0090-5364}, mr={0348906}}
\end{barticle}
\bptok{imsref}%
\endbibitem

\bibitem{birbishol2010}
\begin{barticle}[mr]
\bauthor{\bsnm{Birke},~\bfnm{Melanie}\binits{M.}},
\bauthor{\bsnm{Bissantz},~\bfnm{Nicolai}\binits{N.}} \AND
\bauthor{\bsnm{Holzmann},~\bfnm{Hajo}\binits{H.}}
(\byear{2010}).
\btitle{Confidence bands for inverse regression models}.
\bjournal{Inverse Problems}
\bvolume{26}
\bpages{115020, 18}.
\bid{doi={10.1088/0266-5611/26/11/115020}, issn={0266-5611}, mr={2732912}}
\end{barticle}
\bptok{imsref}%
\endbibitem

\bibitem{bisbir2009}
\begin{barticle}[mr]
\bauthor{\bsnm{Bissantz},~\bfnm{Nicolai}\binits{N.}} \AND
\bauthor{\bsnm{Birke},~\bfnm{Melanie}\binits{M.}}
(\byear{2009}).
\btitle{Asymptotic normality and confidence intervals for inverse regression models with convolution-type operators}.
\bjournal{J. Multivariate Anal.}
\bvolume{100}
\bpages{2364--2375}.
\bid{doi={10.1016/j.jmva.2009.04.004}, issn={0047-259X}, mr={2560377}}
\end{barticle}
\bptok{imsref}%
\endbibitem

\bibitem{bisdumholmun2007}
\begin{barticle}[mr]
\bauthor{\bsnm{Bissantz},~\bfnm{Nicolai}\binits{N.}},
\bauthor{\bsnm{D{\"u}mbgen},~\bfnm{Lutz}\binits{L.}},
\bauthor{\bsnm{Holzmann},~\bfnm{Hajo}\binits{H.}} \AND
\bauthor{\bsnm{Munk},~\bfnm{Axel}\binits{A.}}
(\byear{2007}).
\btitle{Non-parametric confidence bands in deconvolution density estimation}.
\bjournal{J. R. Stat. Soc. Ser. B Stat. Methodol.}
\bvolume{69}
\bpages{483--506}.
\bid{doi={10.1111/j.1467-9868.2007.599.x}, issn={1369-7412}, mr={2323764}}
\end{barticle}
\bptok{imsref}%
\endbibitem

\bibitem{bishohmunruy2007}
\begin{barticle}[mr]
\bauthor{\bsnm{Bissantz},~\bfnm{N.}\binits{N.}},
\bauthor{\bsnm{Hohage},~\bfnm{T.}\binits{T.}},
\bauthor{\bsnm{Munk},~\bfnm{A.}\binits{A.}} \AND
\bauthor{\bsnm{Ruymgaart},~\bfnm{F.}\binits{F.}}
(\byear{2007}).
\btitle{Convergence rates of general regularization methods for statistical inverse problems and applications}.
\bjournal{SIAM J. Numer. Anal.}
\bvolume{45}
\bpages{2610--2636}.
\bid{doi={10.1137/060651884}, issn={0036-1429}, mr={2361904}}
\end{barticle}
\bptok{imsref}%
\endbibitem

\bibitem{cavalier2000}
\begin{barticle}[mr]
\bauthor{\bsnm{Cavalier},~\bfnm{Laurent}\binits{L.}}
(\byear{2000}).
\btitle{Efficient estimation of a density in a problem of tomography}.
\bjournal{Ann. Statist.}
\bvolume{28}
\bpages{630--647}.
\bid{doi={10.1214/aos/1016218233}, issn={0090-5364}, mr={1790012}}
\end{barticle}
\bptok{imsref}%
\endbibitem

\bibitem{cavalier2008}
\begin{barticle}[mr]
\bauthor{\bsnm{Cavalier},~\bfnm{L.}\binits{L.}}
(\byear{2008}).
\btitle{Nonparametric statistical inverse problems}.
\bjournal{Inverse Problems}
\bvolume{24}
\bpages{034004, 19}.
\bid{doi={10.1088/0266-5611/24/3/034004}, issn={0266-5611}, mr={2421941}}
\end{barticle}
\bptok{imsref}%
\endbibitem

\bibitem{cavtsy2002}
\begin{barticle}[mr]
\bauthor{\bsnm{Cavalier},~\bfnm{Laurent}\binits{L.}} \AND
\bauthor{\bsnm{Tsybakov},~\bfnm{Alexandre}\binits{A.}}
(\byear{2002}).
\btitle{Sharp adaptation for inverse problems with random noise}.
\bjournal{Probab. Theory Related Fields}
\bvolume{123}
\bpages{323--354}.
\bid{doi={10.1007/s004400100169}, issn={0178-8051}, mr={1918537}}
\end{barticle}
\bptok{imsref}%
\endbibitem

\bibitem{clakei2003}
\begin{barticle}[mr]
\bauthor{\bsnm{Claeskens},~\bfnm{Gerda}\binits{G.}} \AND
\bauthor{\bparticle{Van} \bsnm{Keilegom},~\bfnm{Ingrid}\binits{I.}}
(\byear{2003}).
\btitle{Bootstrap confidence bands for regression curves and their derivatives}.
\bjournal{Ann. Statist.}
\bvolume{31}
\bpages{1852--1884}.
\bid{doi={10.1214/aos/1074290329}, issn={0090-5364}, mr={2036392}}
\end{barticle}
\bptok{imsref}%
\endbibitem

\bibitem{enghanneu1996}
\begin{bbook}[mr]
\bauthor{\bsnm{Engl},~\bfnm{Heinz~W.}\binits{H.W.}},
\bauthor{\bsnm{Hanke},~\bfnm{Martin}\binits{M.}} \AND
\bauthor{\bsnm{Neubauer},~\bfnm{Andreas}\binits{A.}}
(\byear{1996}).
\btitle{Regularization of Inverse Problems}.
\bseries{Mathematics and Its Applications}
\bvolume{375}.
\blocation{Dordrecht}:
\bpublisher{Kluwer Academic}.
\bid{doi={10.1007/978-94-009-1740-8}, mr={1408680}}
\end{bbook}
\bptok{imsref}%
\endbibitem

\bibitem{eubspe1993}
\begin{barticle}[mr]
\bauthor{\bsnm{Eubank},~\bfnm{R.~L.}\binits{R.L.}} \AND
\bauthor{\bsnm{Speckman},~\bfnm{P.~L.}\binits{P.L.}}
(\byear{1993}).
\btitle{Confidence bands in nonparametric regression}.
\bjournal{J. Amer. Statist. Assoc.}
\bvolume{88}
\bpages{1287--1301}.
\bid{issn={0162-1459}, mr={1245362}}
\end{barticle}
\bptok{imsref}%
\endbibitem

\bibitem{folland1984}
\begin{bbook}[mr]
\bauthor{\bsnm{Folland},~\bfnm{Gerald~B.}\binits{G.B.}}
(\byear{1984}).
\btitle{Real Analysis: Modern Techniques and Their Applications}.
\bseries{Pure and Applied Mathematics (New York)}.
\blocation{New York}:
\bpublisher{Wiley}.
\bid{mr={0767633}}
\end{bbook}
\bptok{imsref}%
\endbibitem

\bibitem{ginnic2010}
\begin{barticle}[mr]
\bauthor{\bsnm{Gin{\'e}},~\bfnm{Evarist}\binits{E.}} \AND
\bauthor{\bsnm{Nickl},~\bfnm{Richard}\binits{R.}}
(\byear{2010}).
\btitle{Confidence bands in density estimation}.
\bjournal{Ann. Statist.}
\bvolume{38}
\bpages{1122--1170}.
\bid{doi={10.1214/09-AOS738}, issn={0090-5364}, mr={2604707}}
\end{barticle}
\bptok{imsref}%
\endbibitem

\bibitem{hall1992b}
\begin{barticle}[mr]
\bauthor{\bsnm{Hall},~\bfnm{Peter}\binits{P.}}
(\byear{1992}).
\btitle{On bootstrap confidence intervals in nonparametric regression}.
\bjournal{Ann. Statist.}
\bvolume{20}
\bpages{695--711}.
\bid{doi={10.1214/aos/1176348652}, issn={0090-5364}, mr={1165588}}
\end{barticle}
\bptok{imsref}%
\endbibitem

\bibitem{hall1993}
\begin{barticle}[auto:STB|2014/01/06|10:16:28]
\bauthor{\bsnm{Hall},~\bfnm{P.}\binits{P.}}
(\byear{1993}).
\btitle{On Edgeworth expansion and bootstrap confidence bands in nonparametric regression}.
\bjournal{J. R. Stat. Soc. Ser. B Stat. Methodol.}
\bvolume{55}
\bpages{291--304}.
\end{barticle}
\bptok{imsref}%
\endbibitem

\bibitem{kaisom2005}
\begin{bbook}[mr]
\bauthor{\bsnm{Kaipio},~\bfnm{Jari}\binits{J.}} \AND
\bauthor{\bsnm{Somersalo},~\bfnm{Erkki}\binits{E.}}
(\byear{2010}).
\btitle{Statistical and Computational Inverse Problems}.
\blocation{Berlin}:
\bpublisher{Springer}.
\bptnote{check year}%
\end{bbook}
\bptok{imsref}%
\endbibitem

\bibitem{khosh2002}
\begin{bbook}[mr]
\bauthor{\bsnm{Khoshnevisan},~\bfnm{Davar}\binits{D.}}
(\byear{2002}).
\btitle{Multiparameter Processes}.
\bseries{Springer Monographs in Mathematics}.
\blocation{New York}:
\bpublisher{Springer}.
\bnote{An introduction to random fields}.
\bid{mr={1914748}}
\end{bbook}
\bptok{imsref}%
\endbibitem

\bibitem{Lounic2011}
\begin{barticle}[mr]
\bauthor{\bsnm{Lounici},~\bfnm{Karim}\binits{K.}} \AND
\bauthor{\bsnm{Nickl},~\bfnm{Richard}\binits{R.}}
(\byear{2011}).
\btitle{Global uniform risk bounds for wavelet deconvolution estimators}.
\bjournal{Ann. Statist.}
\bvolume{39}
\bpages{201--231}.
\bid{doi={10.1214/10-AOS836}, issn={0090-5364}, mr={2797844}}
\end{barticle}
\bptok{imsref}%
\endbibitem

\bibitem{mairuy1996}
\begin{barticle}[mr]
\bauthor{\bsnm{Mair},~\bfnm{Bernard~A.}\binits{B.A.}} \AND
\bauthor{\bsnm{Ruymgaart},~\bfnm{Frits~H.}\binits{F.H.}}
(\byear{1996}).
\btitle{Statistical inverse estimation in {H}ilbert scales}.
\bjournal{SIAM J. Appl. Math.}
\bvolume{56}
\bpages{1424--1444}.
\bid{doi={10.1137/S0036139994264476}, issn={0036-1399}, mr={1409127}}
\end{barticle}
\bptok{imsref}%
\endbibitem

\bibitem{neumann1998}
\begin{barticle}[mr]
\bauthor{\bsnm{Neumann},~\bfnm{Michael~H.}\binits{M.H.}}
(\byear{1998}).
\btitle{Strong approximation of density estimators from weakly dependent observations by density estimators from independent observations}.
\bjournal{Ann. Statist.}
\bvolume{26}
\bpages{2014--2048}.
\bid{doi={10.1214/aos/1024691367}, issn={0090-5364}, mr={1673288}}
\end{barticle}
\bptok{imsref}%
\endbibitem

\bibitem{neupol1998}
\begin{barticle}[mr]
\bauthor{\bsnm{Neumann},~\bfnm{Michael~H.}\binits{M.H.}} \AND
\bauthor{\bsnm{Polzehl},~\bfnm{J{\"o}rg}\binits{J.}}
(\byear{1998}).
\btitle{Simultaneous bootstrap confidence bands in nonparametric regression}.
\bjournal{J. Nonparametr. Statist.}
\bvolume{9}
\bpages{307--333}.
\bid{doi={10.1080/10485259808832748}, issn={1048-5252}, mr={1646905}}
\end{barticle}
\bptok{imsref}%
\endbibitem

\bibitem{owen2005}
\begin{bincollection}[mr]
\bauthor{\bsnm{Owen},~\bfnm{Art~B.}\binits{A.B.}}
(\byear{2005}).
\btitle{Multidimensional variation for quasi-{M}onte {C}arlo}.
In \bbooktitle{Contemporary Multivariate Analysis and Design of Experiments}.
\bseries{Ser. Biostat.}
\bvolume{2}
\bpages{49--74}.
\blocation{Hackensack, NJ}: \bpublisher{World Sci. Publ}.
\bid{mr={2271076}}
\end{bincollection}
\bptok{imsref}%
\endbibitem

\bibitem{parpar1973}
\begin{barticle}[mr]
\bauthor{\bsnm{Paranjape},~\bfnm{S.~R.}\binits{S.R.}} \AND
\bauthor{\bsnm{Park},~\bfnm{C.}\binits{C.}}
(\byear{1973}).
\btitle{Laws of iterated logarithm of multiparameter {W}iener processes}.
\bjournal{J. Multivariate Anal.}
\bvolume{3}
\bpages{132--136}.
\bid{issn={0047-259X}, mr={0326852}}
\end{barticle}
\bptok{imsref}%
\endbibitem

\bibitem{piterbarg1996}
\begin{bbook}[auto:STB|2014/01/06|10:16:28]
\bauthor{\bsnm{Piterbarg},~\bfnm{V.~I.}\binits{V.I.}}
(\byear{1996}).
\btitle{Asymptotic Methods in the Theory of Gaussian Processes and Fields}.
\bseries{Translations of Mathematical Monographs}
\bvolume{148}.
\blocation{Providence, RI}:
\bpublisher{Amer. Math. Soc.}
\bnote{Translated from the Russian by V.V. Piterbarg, Revised by the author}.
\end{bbook}
\bptok{imsref}%
\endbibitem

\bibitem{rio1993}
\begin{barticle}[mr]
\bauthor{\bsnm{Rio},~\bfnm{Emmanuel}\binits{E.}}
(\byear{1993}).
\btitle{Strong approximation for set-indexed partial sum processes via {KMT} constructions. {I}}.
\bjournal{Ann. Probab.}
\bvolume{21}
\bpages{759--790}.
\bid{issn={0091-1798}, mr={1217564}}
\end{barticle}
\bptok{imsref}%
\endbibitem

\bibitem{saitoh1997}
\begin{bbook}[mr]
\bauthor{\bsnm{Saitoh},~\bfnm{S.}\binits{S.}}
(\byear{1997}).
\btitle{Integral Transforms, Reproducing Kernels and Their Applications}.
\bseries{Pitman Research Notes in Mathematics Series}
\bvolume{369}.
\blocation{Harlow}:
\bpublisher{Longman}.
\bid{mr={1478165}}
\end{bbook}
\bptok{imsref}%
\endbibitem

\bibitem{smirnov1950}
\begin{barticle}[mr]
\bauthor{\bsnm{Smirnov},~\bfnm{N.~V.}\binits{N.V.}}
(\byear{1950}).
\btitle{On the construction of confidence regions for the density of distribution of random variables}.
\bjournal{Doklady Akad. Nauk SSSR (N.S.)}
\bvolume{74}
\bpages{189--191}.
\bid{mr={0037494}}
\end{barticle}
\bptok{imsref}%
\endbibitem

\bibitem{truj}
\begin{barticle}[auto:STB|2014/01/06|10:16:28]
\bauthor{\bsnm{Trujillo},~\bfnm{I.}\binits{I.}},
\bauthor{\bsnm{Erwin},~\bfnm{P.}\binits{P.}},
\bauthor{\bsnm{Ramos},~\bfnm{A.~A.}\binits{A.A.}} \AND
\bauthor{\bsnm{Graham},~\bfnm{A.~W.}\binits{A.W.}}
(\byear{2004}).
\btitle{Evidence for a new elliptical-galaxy paradigm: S\'ersic and core galaxies}.
\bjournal{Astron. J.}
\bvolume{127}
\bpages{1917--1942}.
\end{barticle}
\bptok{imsref}%
\endbibitem

\bibitem{xia1998}
\begin{barticle}[mr]
\bauthor{\bsnm{Xia},~\bfnm{Yingcun}\binits{Y.}}
(\byear{1998}).
\btitle{Bias-corrected confidence bands in nonparametric regression}.
\bjournal{J. R. Stat. Soc. Ser. B Stat. Methodol.}
\bvolume{60}
\bpages{797--811}.
\bid{doi={10.1111/1467-9868.00155}, issn={1369-7412}, mr={1649488}}
\end{barticle}
\bptok{imsref}%
\endbibitem

\end{thebibliography}
\end{document}